%% file: Injectives_obstruct_Fourier-Mukai_functors.tex
\documentclass[11pt,letterpaper]{amsart}
\usepackage[utf8]{inputenc}
\usepackage{lmodern}

\usepackage{extdash}
\usepackage{amsmath,amsthm,amsfonts,amssymb, microtype}
\usepackage{xr-hyper}
\usepackage[colorlinks=true]{hyperref}
\usepackage{verbatim}
\usepackage{pdfsync}
\usepackage{stmaryrd}
\usepackage{marvosym}
\usepackage{colonequals}
\usepackage{mathrsfs}

\usepackage{mathtools}
\usepackage{enumitem} 
\usepackage{tikz}
\usetikzlibrary{arrows,calc,positioning}

\usetikzlibrary{backgrounds}

\usepackage{extdash}
\usepackage{amsthm,amsfonts,amssymb,amsmath, microtype}
\usepackage{xr-hyper}
\usepackage[colorlinks=true]{hyperref}
\usepackage{verbatim}
\usepackage{pdfsync}
\usepackage{comment} 

\usetikzlibrary{circuits.logic.IEC}
\usetikzlibrary{cd}
\usepackage{xcolor}
\usetikzlibrary{shapes}
\usetikzlibrary{decorations.markings}
\input{tikzstyledefs.tex}

\newcommand{\I}{\mathrm{I}}

\newcommand{\obj}{\mathrm{obj}}

\newcommand{\id}{\mathrm{Id}}
\newcommand{\im}{\mathrm{im}}
\newcommand{\op}{\mathrm{op}}

\renewcommand{\H}{\mathrm{H}}
\newcommand{\Ext}{\mathrm{Ext}}
\newcommand{\can}{\mathrm{can}}

\newcommand{\Qcoh}{\mathrm{Qcoh}}
\newcommand{\Coh}{\mathrm{coh}}
\newcommand{\modules}{\mathrm{-mod}}
\newcommand{\Rhom}{\mathrm{RHom}}

\newcommand{\A}{\mathcal{A}}
\newcommand{\C}{\mathcal{C}}
\newcommand{\OO}{\mathcal{O}}

\newcommand{\m}{\mathrm{m}}

\newcommand{\G}{\mathcal{G}}

\newcommand{\F}{\mathcal{F}}

\newcommand{\D}{\mathcal{D}}
\newcommand{\X}{\mathcal{X}}
\newcommand{\Y}{\mathcal{Y}}
\newcommand{\M}{\mathcal{M}}

\renewcommand{\I}{\mathcal{I}}

\newcommand{\HH}{\mathrm{HH}}
\newcommand{\PP}{\mathbb{P}}
\newcommand{\B}{\mathcal{B}}

\newcommand{\Ho}{\mathrm{Ho}}

\newcommand{\End}{\mathrm{End}}

\newcommand{\Hom}{\mathrm{Hom}}

\newcommand{\Inj}{\mathrm{Inj}}
\newcommand{\Cat}{\mathrm{Cat}}

\newcommand{\kk}{\Bbbk}

\newcommand{\Ainfty}{\A_\infty}

\newcommand{\dg}{\mathrm{dg}}

\newcommand{\nocontentsline}[3]{}
\newcommand{\tocless}[2]{\bgroup\let\addcontentsline=\nocontentsline#1{#2}\egroup}

\newtheorem{theorem}{Theorem}[section]
\newtheorem{proposition}[theorem]{Proposition}
\newtheorem{lemma}[theorem]{Lemma}
\newtheorem{corollary}[theorem]{Corollary}
\newtheorem*{theorem*}{Theorem}

\theoremstyle{definition}
\newtheorem{definition}[theorem]{Definition}
\newtheorem{lemma/definition}[theorem]{Definition/Lemma}
\newtheorem{example}[theorem]{Example}

\newtheorem{remark}[theorem]{Remark}

\newtheorem{construction}[theorem]{Construction} 
\newtheorem{conjecture}[theorem]{Conjecture}

\begin{document}

\title{Injectives obstruct Fourier-Mukai functors}
\author{Felix K\"ung}
\begin{abstract}
We use injectives as a big tilting object to obstruct liftability of exact functors to the $\dg$-level.

We use the inclusion of injectives into the canonical heart as a replacement for tilting objects in computations of the characteristic morphism. Then we apply this construction to proofs of non-liftability of candidate non-Fourier-Mukai functors, i.e.\ functors that do not admit an $\Ainfty$/$\dg$-lift. This approach allows explicit computation of the obstruction against an $\Ainfty$-lift. We in particular observe that this computation gives for smooth degree $d>2$ hypersurfaces an abundance of non-Fourier-Mukai functors.
\end{abstract} 

\maketitle

\let\thefootnote\relax\footnotetext{The author is a Postdoctoral Research Fellow at the Universit\'e Libre de Bruxelles. The research in this work has been supported by the MIS (MIS/BEJ - F.4545.21) Grant from the FNRS, a mobility grant by the FNRS and an ACR Grant from the Universit\'e Libre de Bruxelles.}

\section{Introduction}

The idea behind Fourier--Mukai functors as categorified correspondences was introduced by S. Mukai \cite{Mukai1981} to understand the relation between the derived category of an abelian variety $\A$ and its dual $\A^\vee$. The classical definition is:

\begin{definition}\cite{HuybrechtsFm}
A functor $$f: \D^b\left( \Coh X\right)\to \D^b\left( \Coh Y\right),$$ respectively $$f\colon \D^b\left(\Qcoh X\right)\to \D^b\left(\Qcoh Y\right),$$ for two schemes $X$ and $Y$ is Fourier-Mukai if there exists a sheaf $M$ in $\D^b\left(\Coh \left(X\times Y\right)\right)$, respectively in $\D^b\left(\Qcoh\left(X\times Y\right)\right)$, such that
$$f\cong\pi_{1,*}\left(M\otimes \pi_2^* \_\right).$$
\end{definition}

 This concept allows for the study of functors from a geometric point of view and in particular the application of geometric arguments and constructions to them. With time it was realized that in fact every known functor between derived categories of smooth projective schemes is Fourier-Mukai. For example, push-forward and pull-back along a morphism of schemes $f$ is induced by taking $M$ to be the structure sheaf of its graph $\OO_{\Gamma_f}$.

 A celebrated result by D.\ Orlov gave a criterion to check whether a functor is Fourier-Mukai:

\begin{theorem}\cite{Orlov1997}
Let $f:\D^b\left(\Coh X\right) \to \D^b\left(\Coh Y\right)$ be a fullyfaithful functor between derived categories of smooth projective schemes. Then $f$ is Fourier-Mukai.
\end{theorem}

Whereas B.\ To\"en realized a beautiful conceptual categorization of the property.

\begin{theorem}\cite{Toen2004}
Let $f:\D^b\left(\Qcoh X\right) \to \D^b\left(\Qcoh Y\right)$ be an exact functor between derived categories of smooth projective schemes. Then $f$ is Fourier-Mukai if and only if there exists a $\dg$-lift of $f$ i.e.\
$$\widehat{f}:\D_{\dg}^b\left(\Qcoh X\right) \to \D_{\dg}^b\left(\Qcoh Y\right),$$
such that $\H^0\left(\widehat{f}\right)\cong f$, where $\D_{\dg}^b\left(\_\right)$ is a $\dg$-enhancement of $\D^b\left(\_\right)$.
\end{theorem}

This theorem in particular allows the interpretation of Fourier-Mukai functors as the geometric analog to bimodules and their represented functors. These observations led to a conjecture that all exact functors between derived categories of smooth projective schemes are Fourier-Mukai.

\begin{conjecture}\cite{Lunts2007}
Let $f: \D^b\left(\Coh X\right)\to \D^b\left(\Coh Y\right)$ be a functor between derived categories of smooth projective schemes. Then $f$ is Fourier-Mukai.
\end{conjecture}

More than ten years later A.\ Rizzardo, M.\ Van den Bergh and A.\ Neeman constructed the first non-Fourier-Mukai functor.

\begin{theorem}\cite{Rizzardo2019}
There exists an exact functor $$\psi_\eta: \D^b\left(\Coh Q_3\right) \to \D^b\left(\Coh \PP^4\right)$$
that is not Fourier-Mukai, where $Q_3 \hookrightarrow \PP^4$ is a smooth three dimensional quadric.
\end{theorem} 

We were able to generalize their argument to higher dimensions in \cite{Kueng2022} using twisted Hodge diamonds, while T.\ Raedschaelders, A.\ Rizzardo and M.\ Van den Bergh \cite{Raedschelders2019} were able to construct similar non-Fourier-Mukai functors by losing control over the target space.

 A big drawback of all these approaches is that they heavily rely on the existence of tilting objects in the source category in order to control liftability, which tend to be rare. We avoid this issue here by replacing the tilting bundle with the subcategory of injectives in the canonical heart of the derived category of quasi-coherent sheaves. This observation allows the computation of the equivariant characteristic morphism with respect to the category of injectives and explicit computations of the induced $\Ainfty$-obstructions. In particular we are able to prove the following very general result:

\setcounter{section}{8}
\setcounter{theorem}{10}

\begin{theorem}
Let $f:X\hookrightarrow Y$ be a divisor and let $0\neq \eta\in \HH^*\left(X,M\right)$ be a Hochschild cohomology class such that $f_*\eta =0$ for $M$ a line bundle on $X$. Then there exists a non-Fourier-Mukai functor
$$\psi_\eta: \D^b\left(\Qcoh X\right) \to \D^b\left(\Qcoh Y\right).$$
\end{theorem}

Explicit computations for smooth hypersurfaces in projective space then yield:

\setcounter{theorem}{12}
\begin{corollary}
Let $f:X\hookrightarrow \PP^{n+1}$ be a smooth degree $d$-hypersurface. Then the space of non-Fourier-Mukai functors arising as deformations of the shape $\psi_\eta$ can be computed explicitly.
\end{corollary}

Finally we show a possible avenue for generalizing our result to coherent sheaves by using extension results by J. Lurie \cite{Lurie} respectively B. To\"en \cite{Toen2004}.

\setcounter{section}{8}
\setcounter{theorem}{11}
\begin{theorem}
Assume that $\psi_\eta$ is continuous. Then the restriction $$\psi_\eta: \D^b\left(\Coh X\right)\to \D^b\left(\Coh Y\right)$$ is not Fourier-Mukai.
\end{theorem} 

\setcounter{section}{1}

\subsection{Proof strategy}

By classical methods $\Ainfty$-categories are equivalent to $\dg$-categories and better behaved with respect to liftability. So we will work with $\Ainfty$-structures in order to obstruct $\dg$/$\Ainfty$-lifts.

Let us briefly sketch the above mentioned previous approaches to non-liftability. These assumed the existence of a tilting bundle $T$ with endomorphism algebra $\Gamma$. After applying the candidate non-Fourier-Mukai functor $\psi_\eta$ this turns $\psi_\eta \left(T\right)$ into a $\Gamma$-equivariant object in $\D^b\left(Y\right)$, whose $\Ainfty$-liftability then can be obstructed.

The key insight in our approach is that the non-liftability of $\psi_\eta \left(T\right)$ boils down to the the functor $\psi_\eta \circ \iota_T : *_\Gamma \to \D^b\left(Y\right)$ not admitting an $\Ainfty$-lift, here we denoted by $*_\Gamma$ the $\kk$-linear category with one object $*$ and endomorphisms $\Gamma$. In particular, since $T$ was assumed to be tilting, the functor $$\iota_T : *_\Gamma \to \D^b\left(X\right)$$ determines functors on $\D^b\left(X\right)$ uniquely. 

Our approach consists of replacing the functor $\iota_T:*_\Gamma \hookrightarrow \D^b\left(X\right)$ with the functor given by the inclusion of the full subcategory of injective objects into the heart $$\iota: \Inj\hookrightarrow \D^b\left(\Qcoh X\right).$$ By classical homological algebra this functor completely controls $\D^b\left(\Qcoh X\right)$ \cite{Weibel1994}.

In order to explain our proof strategy, consider the following diagram

$$\tikz[heighttwo,xscale=4,yscale=2,baseline]{
\node (Inj) at (0,1) {$\Inj$};
\node (DQX) at (1,0) {$\D^b\left(\Qcoh X\right)$};
\node (DXe) at (2,0) {$\D^b\left(\X_\eta\right)$};
\node (DQY) at (3,0) {$\D^b\left(\Qcoh Y\right),$};
\node (DQXdg) at (1,1) {$\D_\dg^b\left(\Qcoh X\right)$};
\node (DXedg) at (2,1) {$\D_\dg^b\left(\X_\eta\right)$};
\node (DQYdg) at (3,1) {$\D_\dg^b\left(\Qcoh Y\right)$.};

\draw[->]
(Inj) edge node[above] {$\iota$} (DQXdg)
(Inj) edge node [below left] {$\iota$} (DQX)
(DQX) edge node[above] {$L$} (DXe)
(DXe) edge node[above] {$\widetilde{f}_*$} (DQY)
(DXedg) edge node[above] {$\widetilde{f}_*$} (DQYdg);
\draw[->, dashed]
(DQXdg) edge (DQX)
(DXedg) edge (DXe)
(DQYdg) edge (DQY);
\draw[->, dashed,bend left = 60]
(Inj) edge node[red] {\begin{Large}\Lightning\end{Large}} (DXedg)
(Inj) edge (DQYdg);
}
$$
where to top row indicates the $\dg$-enhancements of the bottom row. The functor $\psi_\eta$ will be constructed as the composition of $L$ and $\widetilde{f}_*$.
 
As $\iota$ uniquely determines $\dg$-functors from the category $\D^b\left(\Qcoh\left(X\right)\right)$ we can compute its equivariant characteristic morphism in order to obstruct an $\Ainfty$-lift of $L$. This yields a non-vanishing $\Ainfty$-obstruction $o_m\left(L\circ \iota\right)$. The plan is to push this obstruction forward. In contrast to \cite{Rizzardo2019} and \cite{Kueng2022} this needs to be done with a lot of care as we do not have access to Serre duality in our setting. In particular we study the concrete construction of $\widetilde{f}_*$ to compute that the obstruction $o_m\left(\psi_\eta \circ \iota \right)$ indeed is independent of choices and coincides with the pushed forward obstruction $f_*o_m\left(L\circ \iota\right)$. This means that there cannot be an $\Ainfty$-lift of $\psi_\eta \circ \iota$. In particular, as $\iota$ is the inclusion of injectives, it is itself a $\dg$-functor and hence we conclude that $\psi_\eta$ cannot have an $\Ainfty$-lift as desired.

\subsection{Structure of the paper}

In Section~\ref{Section Preliminaries Hochschild cohomology} we recall preliminaries on Hochschild cohomology and the characteristic morphism in the classical setting. In particular we recall in Definition~\ref{Definition infinitesimal deformation along Hochschild cohomology class} the construction of an infinitesimal deformation along a Hochschild cohomology class, which is a key player in the construction of $\psi_\eta.$

Then we recall in Section~\ref{Section equivariant objects over categories} how one can interpret equivariant objects as functors out of one object categories. This will be used later to make our arguments on the characteristic morphism of subcategories clearer.

Using these observations, together with insights from \cite{Toen2004}, in Section~\ref{Section Equivariant characteristic morphism for injectives} we compute that the $\Inj$-equivariant characteristic morphism of the subcategory $\Inj$ does not vanish. This construction will be our replacement for the assumption on the existence of a tilting object.

We recall in Section~\ref{Section Construction of Psi} the construction from \cite{Rizzardo2019} of candidate non-Fourier-Mukai functors. In Section~\ref{Section Recollection obstructions} we then collect the properties, definitions and computations of $\Ainfty$-obstructions that we will use later.

Since we are not using tilting objects in the category of coherent sheaves we are not able to use Serre duality as in \cite{Rizzardo2019,Kueng2022} in order to show that the relevant $\Ext$-groups are concentrated in degree $-1$. So we will have to control the $\Ainfty$-structure on our target category more explicitly. For this we study in Section~\ref{Section widetilde f} the explicit structure of the functor $\widetilde{f}$. In particular we completely understand its essential image and its $\Ainfty$-bimodule structure.

We finish in Section~\ref{Section Computing obstructions explicitly} by computing explicitly that the obstructions of $\psi_\eta$ cannot vanish. We then combine this result with computations from \cite{Kueng2022} in order to conclude the existence of a vast amount of non-Fourier-Mukai functors with well-behaved target categories and showing possible arguments for the case of coherent sheaves.

\section*{Acknowledgements}

I would like to thank the SL-Math institute for hosting the semester on ``non-commutative algebraic geometry'' where a lot of this work has been completed. Furthermore I want to show my gratitude to M. Van den Bergh and A. Rizzardo for helpful discussions and deep insights  and \v{S}. \v{S}penko for her awesome support and helpful discussions. I would like to thank B. Briggs and A. C\v{a}ld\v{a}raru for their insights on Hochschild cohomology, $\Ainfty$-structures and formality of sheaves over source spaces. Finally I are deeply grateful for the insight of J. Holstein on continuation of functors allowing the possible restriction to coherent sheaves. 

\section*{Notation}

Let $\kk$ denote a field of characteristic zero. We consider all schemes and categories to be over $\kk$. For a $\kk$-linear category $\X$, denote by $\D^{b}\left(\X\right)$ the bounded derived category of left modules over $\X$. For $X$ a scheme we denote by $\D^b\left(X\right)$ the bounded derived category of coherent sheaves. Denote by $\D_\dg^b\left(\_\right)$ the (canonical) dg-enhancement of $\D^b\left(\_\right)$. We denote by $\D^b_\C\left(\_\right)$ the full subcategory of $\D^b\left(\_\right)$ with cohomology objects in $\C$. We also use the slightly strange convention of $f_*$ for the pullback of $\kk$-linear categories in order to mimic the geometric setting and improve readability. We will also refer to $\dg$ respectively $\Ainfty$-functors between smooth geometric categories as Fourier Mukai functors following the insight of B. To\"en \cite{Toen2004}. Finally we use the shorthand 

\begin{align*}
\left[f_i,m_k\right]:=& \sum_{\mathclap{r+t=i-1}}\left(-1\right)^{r+kt}f_i\left(\id^{\otimes r}\otimes \m_k \otimes \id^{\otimes t}\right)\\
&-\sum_{\mathclap{1\le j\le k}}\left(-1\right)^{ij} \m_k \left(\alpha_1\otimes ... \otimes f_i\left(\alpha_{j+1},...,\alpha_{j+i} \right) \otimes...\otimes \alpha_{k+i} \right)
\end{align*}

for the Gerstenhaber bracket of a morphism $\A^{\otimes i} \to \M$ and $\M$ an $\Ainfty$-$\A$-module.

\section{Preliminaries}

We start by recalling some basic notions needed for our constructions and arguments. For the basics on $\Ainfty$-structures we refer to the notes by B. Keller \cite{KellerAinfty}, the work of H. Lefevre \cite{Lefevre} and P. Seidel \cite{Seidel2008}. For the construction of $\dg$-hulls we refer to \cite{Kapustin2009}.

\subsection{Hochschild cohomology and characteristic morphism}\label{Section Preliminaries Hochschild cohomology}

As we will work a lot with obstructions in Hochschild cohomology we recall the basic constructions and notions here. For a more in depth discussion on the incarnations of Hochschild cohomology we use here, we refer the reader to the original work of G. Hochschild \cite{Hochschild}, R. Buchweitz and H. Flenner \cite{Buchweitz2006} respectively R. Swan \cite{Swan1996} for its geometric incarnation, and for its relation to deformation theory to the work of M. Van den Bergh and W. Lowen \cite{Lowen2005}. Finally we want to refer for a treatment on the characteristic morphism to the work of B. Briggs and V. Gelinas \cite{Briggs2017}.

\begin{definition}
Let $\X$ be a small $\kk$-linear category and $\M$ a $\kk$-central $\X$-bimodule. Then the Hochschild complex of $\X$ with values in $\M$ is given by the complex
$$\mathrm{CC}^n\left(\X,\M\right):=\prod_{\mathclap{X_0,...,X_n \in \obj \left(\X\right)}}\Hom\left(\X\left(X_0,X_1\right)\otimes_\kk \cdots \otimes_\kk \X\left(\X_{n-1},\X_n\right),\M\left(X_0,X_n\right)\right)$$
with differential given by
\begin{align*}
d f\left(x_1\otimes \cdots \otimes x_{n+1}\right):=& x_1f\left(x_2\otimes \cdots \otimes x_{n+1}\right)\\
&+\sum_{i=1}^n\left(-1\right)^i f\left(x_1\otimes\cdots \otimes x_i x_{i+1}\otimes \cdots \otimes x_{n+1}\right)\\
&+\left(-1\right)^{n+1}f\left(x_1\otimes \cdots \otimes x_n\right)x_{n+1}.
\end{align*}
The Hochschild cohomology $\HH^*\left(\X,\M\right)$ is the cohomology of $\mathrm{CC}^\cdot\left(\X,\M\right)$.
\end{definition}

The following classical Lemma~\ref{Lemma Hochschild cohomology is ext in bimodules} allows for a more streamlined construction and will help us interpret the Hochschild cohomology in more generality later.

\begin{lemma}\label{Lemma Hochschild cohomology is ext in bimodules}
Let $\X$ be a $\kk$-linear category and let $\M$ be a $\kk$-central bimodule over $\X$. Then the Hochschild cohomology of $\X$ with values in $\M$ is $$\HH^*\left(\X,\M\right):=\Ext^*_{\X\otimes \X^\op}\left(\X,\M\right).$$ 
\end{lemma}

\begin{definition}\label{Definition tensor product with small K-linear category}
Let $\X$ be a small $\Ainfty$-category, and let $\I$ be a small $\kk$-linear category. Then we define the $\Ainfty$-category $\X\otimes_\kk \I$ to consist of the objects $$\left(x,i\right)\in \obj\left(\X\right)\times \obj\left(\I\right),$$ morphism spaces given by
$$\left(\X\otimes_\kk \I\right)\left(\left(x,i\right)\left(y,j\right)\right):=\X\left(x,y\right)\otimes_\kk \I\left(i,j\right)$$
and higher compositions 
$$\m_{k,\X\otimes\I}\left(\left(x_1\otimes_\kk i_1\right),...,\left(x_k\otimes_\kk i_k\right)\right):=\m_{k,\X}\left(x_1,...,x_k\right)\otimes i_1 \cdots i_k.$$
\end{definition}

\begin{remark}
In bigger generality constructing tensor products of $\Ainfty$-categories needs to be done with a lot of care \cite{Lowen2020}. However, we will only consider tensor products with $\kk$-linear categories, which are better behaved and allow our naive approach to work.
\end{remark}

The restriction in the following definition to $m\geq 3$ is not necessary in general. However, as we will have naturally this assumption and it makes the formulas more concise we present this special case.

\begin{definition}\label{Definition infinitesimal deformation along Hochschild cohomology class}
Let $\X$ be a $\kk$-linear category and let $\M$ be a $\kk$-central $\X$-bimodule.
Let $\eta \in \mathrm{CC}^{m}\left( \X,\M\right)$ be a cocycle with $m\geq 3$. Then the infinitesimal deformation of $\X$ along $\eta$ is the $\Ainfty$-category $\X_\eta$ given by objects $\obj \X_\eta =\obj \X$, morphisms
$$\X_\eta \left(x,y \right):= \X\left(x,y\right)\oplus \M\left(x,y\right)\left[2-m\right]$$
and non-zero multiplications given by 
\begin{align*}
\m_2: \left(\X\left(x,y\right)\oplus \M\left(x,y\right)\left[2-m\right]\right)^{\otimes 2}&\to \X\left(x,y\right)\oplus \M\left(x,y\right)\left[2-m\right]\\
\left(x,m\right)\otimes \left(x',m'\right)&\mapsto \left(x x',x m' + m x'\right)\\
\m_m:\left(\X\left(x,y\right)\oplus \M\left(x,y\right)\left[2-m\right]\right)^{\otimes m} &\to \X\left(x,y\right)\oplus \M\left(x,y\right)\left[2-m\right]\\
x_1\otimes... \otimes x_m &\mapsto \eta\left(x_1,...,x_m\right).
\end{align*}

This construction comes with a canonical morphism that is the identity on $\X$ and zero on $\M$,
\begin{equation}\label{canonical morphism for deformation}
\mathrm{can}:\X_\eta \to \X.
\end{equation}
\end{definition}

\begin{proposition}\cite[Lemma~6.1.1.]{Rizzardo2019}\label{Proposition construction X eta independent of choice of cocycle}
Let $\X$ be a $\kk$-linear category, $\M$ a $\kk$-central $\X$-bimodule and $\eta,\mu \in \mathrm{CC}^{m}\left( \X,\M\right)$ such that $\eta=\mu$ in $\HH^*\left(\X,\M\right)$. Then 
$$\X_\eta\cong \X_\mu .$$
\end{proposition}

\begin{remark}\label{Remark comparison of deformations}
Let $\X$ and $\I$ be $\kk$-linear categories and $\M$ a $\kk$-central bimodule. Then there is a morphism of Hochschild complexes 
\begin{align*}
\mathrm{CC}^*\left(\X,\M\right)&\to \mathrm{CC}^*\left(\X\otimes \I,\M\otimes \I\right)\\
\eta&\mapsto \eta \cup \id 
\end{align*}
One can compute that this morphism is compatible with deformations, i.e.
$$\X_\eta\otimes \I\cong \left(\X\otimes \I\right)_{\eta \cup \id}.$$
\end{remark}

\begin{definition}\label{Definition characteristic morphism naive}
Let $\X$ be a $\kk$-linear category, $\M$ a $\kk$-central bimodule and let $\G$ be a $\X$-module. Then $\_\otimes \G$ induces a  map called the characteristic morphism $$c_{\G}:\HH^*\left(\X,\M\right)\cong \Ext^*_{\X\otimes \X^\op}\left(\X,\M\right) \xrightarrow{\otimes_\X \G} \Ext^*_{\X}\left(\G,\M\otimes \G\right).$$

Analogously we have the dual characteristic morphism
$$c^*_{\G}:\HH^*\left(\X,\M\right)\cong \Ext^*_{\X\otimes \X^\op}\left(\X,\M\right) \xrightarrow{\Rhom_\X\left(\_, \G\right)} \Ext^*_{\X}\left(\Rhom\left(\M,\G\right), \G\right),$$
induced by $\Rhom_X\left(\_,\G\right)$.
\end{definition}

\subsection{Equivariant objects over categories}\label{Section equivariant objects over categories}

In this section we recall the basic notions of equivariant objects and how to interpret them as functors from a small category. We will apply this in the next section to compute the equivariant characteristic morphisms of these objects.

We start by recalling the basic definition.

\begin{definition}\label{definition equivariant object classical}
Let $\C$ be a $\kk$-linear category and let $\Gamma$ be a $\kk$-algebra. A $\Gamma$-equivariant object $X\in \C$ is an object equipped with an action morphism of $\kk$-algebras
$$\varphi :\Gamma \to \C\left(X,X\right).$$

A morphism between $\Gamma$-equivariant objects $X$ and $Y$ is a morphism $f: X \to Y$ in $\C$ such that $$\varphi_Y\left(\gamma\right) \circ f = f\circ \varphi_X \left(\gamma\right).$$
\end{definition}

If there is a $\Gamma$-action on $T$ the characteristic morphism from Definition~\ref{Definition characteristic morphism naive} admits an equivariant incarnation keeping track of the $\Gamma$-action.

\begin{definition}\label{Definition equivariant characteristic morphism naive}
Let $\X$ be a $\kk$-linear category, $\M$ a $\kk$-central bimodule and let $\G$ be a $\Gamma$-equivariant $\X$-module. Then $\_\otimes \G$ induces the $\Gamma$-equivariant characteristic morphism $$c_{\Gamma,\G}:\HH^*\left(\X,\M\right) \to \Ext^*_{\X\otimes \Gamma^{\op}}\left(\G,\M\otimes \G\right).$$
While $\Rhom\left(\_, \G\right)$ induces the $\Gamma$-equivariant dual characteristic morphism $$c^*_{\Gamma,\G}:\HH^*\left(\X,\M\right) \to \Ext^*_{\X\otimes \Gamma^{\op}}\left(\Rhom\left(\M,\G\right), \G\right).$$
\end{definition}

\begin{lemma}\label{Lemma equivariant classical is generalized by categorical}
Let $\C$ be a $\kk$-linear category and let $\Gamma$ be a $\kk$-algebra. Then a $\Gamma$-equivariant object $X$ is equivalently a functor from the one object category $*_\Gamma$ with endomorphisms $\Gamma$, 
\begin{align*}
\iota_{X,\Gamma}:*_\Gamma &\to \C \\
*&\mapsto X	.
\end{align*}
Furthermore this naturally extends to an equivalence of categories between $\Gamma$-equivariant objects in $\C$ and the functor category $\mathrm{Fun}\left(*_\Gamma,\C\right)$.
\begin{proof}
Let $\iota_{X,\Gamma}$ be a functor of the given shape, then it defines a $\Gamma$-equivariant object by taking $X=\iota_{X,\Gamma}\left(*\right)$ with action given by 
$$\iota_{X,\Gamma}: \End \left(*\right)= \Gamma \longrightarrow X = \iota_{X,\Gamma}.$$

Furthermore we can make a $\Gamma$-equivariant object $X$ into a functor via fixing the image of $*$ to be $X$ and the structure morphism given by the map on morphism spaces. This gives us two inverse constructions.

Finally, the naturality requirement for a morphism $\alpha: \iota_{X,\Gamma}\to \iota_{Y,\Gamma}$ in the functor category boils down to
\begin{align*}
\alpha \circ \iota_{X,\Gamma}\left(\gamma \right) = \iota_{Y,\Gamma}\left(\gamma\right)\circ \alpha,
\end{align*}
which is precisely the identity from Definition~\ref{definition equivariant object classical}.
So we may identify the two categories.
\end{proof}
\end{lemma}

In particular we may consider equivariant object over small $\kk$-linear categories by considering functors.

\begin{definition}\label{Definition Equivariant object over category}
Let $\C$ be a $\kk$-linear category and $\I$ a small $\kk$-linear category. An $\I$-equivariant object in $\C$ is a $\kk$-linear functor
$$\I\to \C.$$
We in particular use the following notation for the category of equivariant objects over $\I$ in $\C$ by
$$\C_\I:= \mathrm{Fun}_\kk\left(\I,\C\right).$$
\end{definition}

%\begin{remark}\label{Remark reinterpretation original proof}
%Using the definition of an equivariant object being a functor we can reinterpret the proofs in \cite{Rizzardo2019},\cite{Raedschelders2019} and \cite{Kueng2022} of the non-existence of lifts:

%Consider the composition 
%$$*_\Gamma \hookrightarrow \D^b\left(X\right)\xrightarrow{\psi_\eta} \D^b\left(Y\right).$$
%The above papers conclude the non-existence of a $\dg$-lift by arguing that the $\Gamma$-equivariant object defined by the above functor does not admit an $\Ainfty$-lift to $\Y\otimes \Gamma$. In particular we can reinterpret this as the above composition not being $\dg$. However, as their equivariant object is a tilting object the functor $$*_\Gamma \hookrightarrow \D^b\left(X\right)$$ is $\dg$ and so $\psi_\eta$ cannot be $\dg$.
%\end{remark}

%We will use this more general incarnation of the above strategy which will allow us to consider equivariant objects over the full subcategory of injectives. For this we need to understand the equivariant characteristic morphisms of subcategories better.

\section{Equivariant characteristic morphism of Injectives}\label{Section Equivariant characteristic morphism for injectives}

We study the characteristic morphism introduced in the previous chapter for the subcategory of injectives. We start by discussing how we can reinterpret Hochschild cohomology and the characteristic morphism by the work of \cite{Toen2004}. We then use these insights to prove that we can treat the full subcategory of injectives similarly to a tilting object.  

\subsection{Interpreting characteristic morphisms of a subcategory as pullback}\label{Subsection categorical interpretation of characteristic morphism}

We use observations from \cite{Toen2004} to study the characteristic morphism 
$$c_G:\HH^*\left(\C,\M\right)\to \Ext_\C^*\left(G,\M\otimes_\C G\right)$$
and the equivariant characteristic morphism
$$c_{\Gamma,T}:\HH^*\left(\C,\M\right)\to \Ext_{\C_\Gamma }^*\left(T,\M\otimes_\C T\right).$$

We start by recalling from the previous Section~\ref{Section equivariant objects over categories} that an object $G\in \C$ is equivalently a functor from the one object category $*_\kk$ with endomorphisms $\kk$ into $\C$. In particular we may consider for an object $G\in \C$ the functor
$$\iota_G: *_\kk \to \C, $$
respectively for a $\Gamma$-equivariant $T\in \C$ the functor from the one object category $*_\Gamma$ with endomorphism $\Gamma$ into $\C$. So we may consider in this case $$\iota_T: *_\Gamma \to \C.$$

\begin{definition}
Let $\C$ be an abelian category and $\I\subset \C$ a subcategory. Then denote by
$\iota_I: \I \hookrightarrow \C\subset \D_{\dg}^b\left(\C\right)$ the inclusion of $\I$. 
\end{definition}

\begin{remark}\label{Remark categorical interpretation of characteristic morphism}
By \cite{Toen2004} we can interpret the Hochschlid cohomology of a $\kk$-linear category $\C$ with values in a $\kk$-central bimodule $\M$ as the following $2$-morphism space in the homotopy category of $\dg$-categories.
$$\HH^*\left(\C,\M\right)\cong \Ho\left(\dg-\Cat\right)\left(\D^b_\dg\left(\C\right),\D^b_\dg\left(\C\right)\right)\left(\id,\psi_\M\right)$$
where $\psi_\M$ denotes the $\dg$-functor reperesented by $\M$ i.e.
\begin{align*}
\psi_\M: \D^b\left(\C\right)&\to \D^b\left(\C\right)\\
X&\mapsto \M\otimes X.
\end{align*}

In this interpretation, for a $\Gamma$-equivariant object $T$ the equivariant characteristic morphism turns into application of the pullback functor along the inclusion of the one object category $\iota_T:\Gamma \to \D_\dg^b\left(\Qcoh X\right)$:

$$\tikz[heighttwo,xscale=4,yscale=2,baseline]{
\node (HH) at (0,1) {$\HH^*\left(\C,\M\right)\cong \Ho\left(\dg \Cat\right)\left(\D^b_\dg\left(\C\right),\D^b_\dg\left(\C\right)\right)\left(\id,\psi_\M\right)$};
\node (Ext) at (0,0) {$\Ho\left(\dg \Cat\right)\left(\Gamma,\D^b_\dg\left(\C\right)\right)\left(\iota_T,\psi_\M\circ\iota_\G\right)\cong \Ext_{\C_\Gamma}\left(T,T\otimes \M\right).$};
\draw[->]
(HH) edge node[right] {$c_{\Gamma,T}\cong\iota_{T,*}$} (Ext);
}
$$

In particular we can use this interpretation to compute the equivariant characteristic morphism of the subcategory $\I\subset \C$:
$$\tikz[heighttwo,xscale=4,yscale=2,baseline]{
\node (HH) at (0,1) {$\HH^*\left(\C,\M\right)\cong \Ho\left(\dg\Cat\right)\left(\D^b_\dg\left(\C\right),\D^b_\dg\left(\C\right)\right)\left(\id,\psi_\M\right)$};
\node (Ext) at (0,0) {$\Ho\left(\dg\Cat\right)\left(\I,\D^b_\dg\left(\C\right)\right)\left(\iota,\psi_M\circ\iota\right)\cong \Ext^*_{\C_\I}\left(\iota,\psi_M\circ \iota\right).$};
\draw[->]
(HH) edge node[right] {$c_{\iota_\I,\I}=\iota_{\I,*}$} (Ext);
}
$$
\end{remark}

The above observation is in particular useful as explicit computations for bimodules between $\kk$-linear categories, or, even more, $\Ainfty$/$\dg$-categories can be very unwieldy. However, the above observation allows us to compute the equivariant characteristic morphisms of the $\I$-$\C$ bimodule $\iota_\I$ in special cases.

\subsection{The $\Inj$-equivariant characteristic morphism}

We use the observations on the structure of the characteristic morphism from the last section in order to compute the equivariant characteristic morphism of the inclusion of injective objects. In particular we consider throughout this section the functor 
$$\iota: \Inj \hookrightarrow \C \subset \D_{\dg}\left(\C\right),$$
for $\C$ a category with enough injectives.

\begin{lemma}\label{Lemma inclusion of injectives is dg}
Let $\C$ be a $\kk$-linear category with enough injectives and finite global dimension. Then the inclusion $\iota_*: \Inj \hookrightarrow \C\subset \D^b_{\dg}\left(\C\right)$ is dg.
\begin{proof}
Let $I,J \in \Inj$. Then $\iota \left(I\right)$ and $\iota\left(J\right)$ are injective. In particular we have 
\begin{align*}
\D^b_{\dg}\left(\C\right)\left(\iota I,\iota J\right)&=\D^b_{\dg}\left(\C\right)\left(I,J\right)\\
&=\C\left(I,J\right)
\end{align*}
as graded vector-spaces concentrated in degree $0$. In particular $$\iota: \Inj \hookrightarrow \D^b_\dg\left(\Qcoh\left(X\right)\right)$$ is dg.
\end{proof}
\end{lemma}

\begin{proposition}\label{Proposition inclusion of injectives induces equivalence}
Let $\C$ be a $\kk$-linear category with enough injectives and finite global dimension.
Then the functor $\iota^*: \D^b_\dg\left(\C\right)\to \D^b_\dg\left(\Inj\right)$ is a quasi-equivalence.
\begin{proof}
Since $\C$ has finite global dimension every object in $\D^b\left(\C\right)$ has a finite resolution of injectives. In particular every object $N\in \D^b_\dg\left(\C\right)$ is quasi-isomorphic to a bounded complex of injectives. These complexes are dg-$\Inj$-modules and so $\iota: \Inj\hookrightarrow \C$ induces an equivalence $$\iota_*: \D^b_\dg\left(\C\right)\to \D^b_\dg\left(\Inj\right)$$ as claimed.
\end{proof}
\end{proposition}

\begin{corollary}\label{Corollary characteristic morphism is iso general}
Let $\C$ be a $\kk$-linear category with enough injectives and finite global dimension. Then the equivariant characteristic morphism $$c_{\iota,\Inj}\cong \iota_*: \HH\left(\C,\M\right)\to \Ext_{\C_\Inj}\left(\iota, \psi_\M\circ \iota\right)$$ is an isomorphism.
\begin{proof}
By Section~\ref{Subsection categorical interpretation of characteristic morphism} above we can interpret the equivariant characteristic morphism as the pullback along the inclusion of the full subcategory $$\iota:\Inj \to \D^b_\dg\left(\C\right)$$ in the homotopy category of $\dg$-categories. In particular we need to pick a fibrant replacement of $\Inj$. By \cite[Lemme~2.13]{Tabuada2007} pretriangulated replacements are an incarnation of a fibrant replacement. So we may choose the canonical pretriangulated hull given by $\dg$-modules $\D^b_\dg\left(\Inj\right)$. 
By Proposition~\ref{Proposition inclusion of injectives induces equivalence} $\iota$ induces an equivalence
$$\iota_*: \D^b_\dg\left(\C\right)\to \D^b_\dg\left(\Inj\right)$$
and hence in particular 
$$c_{\iota,\Inj}\cong \iota_*: \HH^*\left(\C,\M\right)\to \Ext^*_{\C_\Inj}\left(\iota, \psi_\M\circ \iota\right)$$
 is an isomorphism.
\end{proof}
\end{corollary}

\begin{corollary}\label{Corollary characteristic morphism is iso}
Let $X$ be a smooth finite dimensional scheme and consider the category $\Qcoh\left(X\right)$. Then the equivariant characteristic morphism $c_{\iota, \Inj}$ is an isomorphism.
\begin{proof}
We have that $\Qcoh\left(X\right)$ has enough injectives. In particular by Corollary~\ref{Corollary characteristic morphism is iso general} we have that 
$$c_{\iota,\Inj}:\HH^*\left(\Qcoh\left(X\right),\M\right) \xrightarrow{\sim}  \Ho\left(\dg\Cat\right)\left(\Inj,\D^b_\dg\left(\Qcoh\left(X\right)\right)\right)\left(\iota,\psi_\M\circ\iota\right)$$ 
is an isomorphism.
\end{proof}
\end{corollary}

\begin{remark}
The above discussion shows that for smooth projective schemes the full subcategory of injectives behaves like a tilting object.
\end{remark}

\section{Construction of $\psi_\eta$}\label{Section Construction of Psi}

We recall the construction of the candidate non-Fourier-Mukai functor $\psi_\eta$. 

In order to apply Definition~\ref{Definition curly X}, Lemma~\ref{Lemma fully faithful embedding into X} and Lemma~\ref{Lemma associated f between k-linear categories} we assume that every quasi-projective scheme $X$ comes equipped with an open affine cover 
$$X=\bigcup_{i=0}^k U_i.$$

We recall the following construction originally introduced by W. Lowen and M. Van den Bergh \cite{Lowen2011}.

\begin{definition}\cite[Definition~4.2]{Lowen2011}\label{Definition curly X}
Let $X=\bigcup_{i=0}^k U_i$ be an open affine cover of a quasiprojective scheme $X$. Consider for  $ \emptyset \neq I\subset \left\{1,...,k\right\}$ the open sets $U_I:=\bigcap_{i\in I}U_i$. Then $\X$ is the category with objects $\mathrm{P}\left\{1,...,k\right\}\setminus \emptyset$, morphisms:
$$\X\left(I,J\right):=\begin{cases}\OO_X\left(U_J\right) & I \subset J\\
0 &\text{else}\end{cases}$$
and composition of morphisms induced by multiplication in the local rings composed with the restriction morphism.
\end{definition}

One can think about $\X\modules$ roughly as the category of presheaves associated to an affine open cover. It comes with the following useful property which is a collection of results from \cite{Rizzardo2019} and which originally appeared in \cite{Lowen2011}.

\begin{lemma}\cite[Lemma~5.2]{Kueng2022} \label{Lemma fully faithful embedding into X}
Let $X$ be smooth quasi-projective. Then there is a fully faithful embedding 
$$w:\D^b\left(\Qcoh X\right)\hookrightarrow \D^b\left(\X\right)$$
with essential image given by $\D^b_{w\Qcoh X}\left(\X\right)$ and a fully faithful embedding
$$W: \Delta_*\D\left(\Qcoh X\right)\hookrightarrow \D\left(\X\otimes \X^{\op}\right),$$
where $\Delta_*\D\left(\Qcoh X\right)$ is the essential image of the direct image of the diagonal $\Delta: X \hookrightarrow X\times X$. In particular we have for quasi-coherent $M$
$$\HH^*\left(X,M\right)\cong \HH^*\left(\X,WM\right).$$
\end{lemma}

\begin{remark}
We will denote $WM$ by $\M$ for readability, similar to $\X$ for the $\kk$-linear category associated to $X$. Observe in particular that we have $W\OO_\Delta \cong \X$ giving rise to the last equality in Lemma~\ref{Lemma fully faithful embedding into X}.
\end{remark}

\begin{corollary}\label{Corollary Inj-equivariant characteristic morphism of a concrete class is non-zero}
Let $\X$ be a $\kk$-linear category associated to a scheme, let $\M=WM$ for $M$ a quasicoherent sheaf and let $0 \neq \eta \in \HH^*\left(\C,\M\right)$. Then we have that its image under the $\Inj$-equivariant characteristic morphism
$c_{\iota,\Inj}\left(\eta\right)$ does not vanish.
\end{corollary}

\begin{lemma}\cite[Lemma~5.4]{Kueng2022}\label{Lemma associated f between k-linear categories}
Let $f: X \hookrightarrow Y$ be a closed embedding of quasi projective schemes. Then we may pull back the open affine cover of $Y$ in the construction of $\Y$ and get an induced functor $\mathfrak{f}: \Y \to \X$ such that the diagram
$$\tikz[heighttwo,xscale=4,yscale=3,baseline]{
\node (DX) at (0,0) {$\D^b\left(\X\right)$};
\node (DQX) at (0,1) {$\D^b\left(\Qcoh X\right)$};
\node (DQY) at (1,1) {$\D^b\left(\Qcoh Y\right)$};
\node (DY) at (1,0) {$\D^b\left(\Y\right),$};

\draw[->]
(DX) edge node[above] {$\mathfrak{f}_*$} (DY)
(DQX) edge node[above left] {$w$} (DX)
(DQY) edge node[above left] {$w$} (DY)
(DQX) edge node[above left] {${f}_*$} (DQY)
;
}
$$
commutes.
\end{lemma}

As we  will mostly be working with the category $\D^b\left(\X\right)$ and $w$ is an exact embedding we will denote the functor $\mathfrak{f}:\X \to \Y$ by $f:\X \to \Y$ for ease of readability.

The following is an application of \cite[Proposition~5.3.1]{Rizzardo2019} to our setting.

\begin{proposition}\cite[Proposition~5.5]{Kueng2022}\label{Proposition construction L}
Let $X$ be smooth projective of dimension $n$ and let $\eta \in \HH^{\geq n+3}\left(X,M\right)$. Then there exists an exact functor
$$L:\D^b\left(\Qcoh X\right)\to \D^b_{wQcoh\left(X\right)}\left(\X_\eta\right).$$
\end{proposition}

We will study the functor arising from the following Lemma~\ref{Lemma construction widetilde f} in more detail in Section~\ref{Section widetilde f}. However, since we want to recall here the basic ingredients to the construction we include its basic construction.

\begin{proposition}\cite[Proposition~8.2.6]{Rizzardo2019}\label{Lemma construction widetilde f}
Let $f: \Y \to \X$ be a $\kk$-linear functor between $\kk$-linear categories $\X$ and $\Y$, $\M$ a $\kk$-central $\X$-bimodule and $\eta\in \HH^{\geq 3}\left(\X,\M\right)$ such that   $0=f_*\eta \in \HH^*\left(\Y,f_*\M\right)$. Then there exists an $\Ainfty$-functor $\widetilde{f}$ making the following diagram commute
$$\tikz[heighttwo,xscale=2.5,yscale=3,baseline]{
\node (X) at (0,0) {$\X$};
\node (Xeta) at (1,1) {$\X_\eta$};
\node (Y) at (2,0) {$\Y$.};

\draw[->]
(Y) edge node[above] {$f$} (X)
(Xeta) edge node[above left] {$\text{can}$} (X)
;
\draw[->,dashed]
(Y) edge node[above right] {$\widetilde{f}$} (Xeta);
}
$$
\end{proposition}

\begin{construction}\cite{Rizzardo2019}
Let $X\hookrightarrow Y$ be the inclusion of a smooth divisor into a smooth projective scheme and let
$$0\neq \eta \in \ker\left(f_*:\HH^m\left(X,M\right)\to \HH^m\left(Y,f_*M\right)\right)$$
with $m\geq \dim X + 3$. 

Then we can construct a functor 

$$\psi_\eta:\D^b\left(\Qcoh X\right)\xrightarrow{L} \D^b_{w\Qcoh X}\left(\X_\eta\right)\xrightarrow{\widetilde{f}_*}\D^b_{w\Qcoh Y}\left(\Y\right)\cong \D^b\left(\Qcoh Y\right),$$
where the functor $L$ is given by Proposition~\ref{Proposition construction L} and $\widetilde{f}_*$ exists by Proposition~\ref{Lemma construction widetilde f}.
\end{construction}

\begin{remark}
Construction~\ref{Lemma construction widetilde f} works as well if $\eta$ is zero, however, the arising functor then trivially admits a lift to the $\Ainfty$-level and so we will not consider this case.
 \end{remark}

\section{Recollection on Obstructions}\label{Section Recollection obstructions}

As we will obstruct $\Ainfty$- respectively $\dg$-lifts we recall the following notions and results from \cite{Rizzardo2019} on colifts along infinitesimal deformations and $\Ainfty$-obstructions:

\begin{definition}\label{Definition lift and colift}
Let $\X$ be a small $\kk$-linear category, assume that $\M$ is a right-flat $\kk$-central bimodule, $\eta \in \HH^{\geq 3}\left(\X,\M\right)$ and let $U \in \X\modules$. 

A lift of $U$ to $\X_\eta$ is an object $V\in \D^b_\dg\left(\X_\eta\right)$, such that there exists an isomorphism in $\X_\eta\modules$:
$$\phi: \H^*\left(V\right)\xrightarrow{\sim} \H^*\left(\X_\eta\right)\otimes_\X U.$$

Dually, if $\M$ is a left-projective $\kk$-central bimodule and $U \in \X\modules$. A colift of $U$ to $\X_\eta$ is an object $V \in \D^b_\dg\left(\X_\eta \right)$, such that there exists an isomorphism in $\X_\eta\modules$:
$$\phi: \H^*\left(V\right)\xrightarrow{\sim}\Hom_\X\left(\H^*\left(\X_\eta\right),U\right).$$
\end{definition}

\begin{theorem}[{\cite[Lemma~8.2.1]{Rizzardo2019}}]\label{Theorem characteristic morphism obstructs lifts and colifts}
The object $U\in \X\modules$ has a lift to $\X_\eta$ if and only if $c_{U}\left(\eta\right)=0$. It has a colift if and only if $c^*_U\left(\eta\right)=0$.
\end{theorem}

For the next Lemma~\ref{Lemma comparison characteristic morphism and dual characteristic morphism} observe that if $M$ is a line-bundle we have that the bimodule $\M:=WM$ is an invertible $\X$-bimodule.

\begin{lemma}[{\cite[Lemma~6.3.1]{Rizzardo2019}}]\label{Lemma comparison characteristic morphism and dual characteristic morphism}
Assume $\M$ is an invertible $\X$-bimodule, then we have that $c_U\left(\eta\right)\neq 0$ if and only if $c^*_U\left(\eta\right)\neq 0$.
\end{lemma}

\begin{lemma}[{\cite[Lemma~8.2.1]{Rizzardo2019}}]\label{Lemma A-infinity obstructions against lift of functors}
Let $f_i: \A \to \mathcal{B}$ be an $\A_i$-functor between two $\Ainfty$-categories. Then there are inductive obstructions 
$$o_{i+1}\left(f_i\right) \in \HH^{i+1}\left(\H^*\left(\A\right),\H^*\left(\mathcal{B}\right)\right)_{1-i}$$
 such that $o_{i+1}\left(f_i\right)$ vanishes if and only if there exists a $\delta$ with $[\delta , m_1]=0$ such that $f_i+\delta$ can be lifted to an $\A_{i+1}$-morphism. 

This obstruction is natural i.e. if we have two $\Ainfty$-functors $G: \A' \to \A$ and $G':\mathcal{B} \to \mathcal{B}'$ we get
$$o_{i+1}\left(G'\circ f_i\circ G\right)=\H^*\left(G'\right)\circ o_{i+1}\left(f_i\right)\circ \H^*\left(G\right).$$
\end{lemma}

The obstructions in Lemma~\ref{Lemma A-infinity obstructions against lift of functors} are given by the $\Ainfty$-functor formulas and the observation that they naturally define Hochschild cohomology classes. In particular they are given by the formulas
\begin{equation}\label{Equation defining obstructions}
o_{i}\left(f\right) = \sum_{k=1}^{i} \pm \left[\m_k,f_{i-k+1}\right] \in \HH^{i+1}\left(\H^*\left(\A\right),\H^*\left(\mathcal{B}\right)\right)_{1-i}
\end{equation}

 \begin{lemma}\label{Theorem Colift and Ainfinity lifts coincide}
Let $\iota \in \D^b \left(\X\otimes \I^{op}\right)$ for a small $\kk$-linear category $\I$. Let $\widehat{\iota}\in \D^b_\mathrm{dg}\left(\X_\eta\otimes \I^{op}\right)$. Then the following are equivalent:
\begin{itemize}
\item The object $\widehat{\iota}$ is a colift of $\iota \in \D^b\left(\X\otimes\I^{op}\right)$.
\item The object $\widehat{\iota}$ is an $\Ainfty$-lift of $L\circ \iota \in \D^b\left(\X_\eta\right)_\I$ i.e. $L\circ \iota: \I \to \X_\eta$. 
\end{itemize}
\begin{proof}
The object $\widehat{\iota}$ is a colift of $\iota$ if and only if we have an isomorphism of $\H^*\left(\X_\eta\right)$-modules:
$$\Phi:\Ext^*_{\X\otimes \I^{op}}\left(\left(\X\otimes \I^{op}\right)_{\eta \cup 1}, \iota\right)\xrightarrow{\sim}\H^*\left(\widehat{\iota}\right).$$
In particular it suffices, since by Remark~\ref{Remark comparison of deformations} $\X_\eta\otimes \I^{op}\cong \left(\X\otimes \I^{op}\right)_{\eta\cup 1}$, to find an isomorphism
$$\Ext^*_{\X\otimes \I^{op}}\left(\left(\X\otimes \I^{op}\right)_{\eta \cup 1}, \iota\right)\xrightarrow{\sim}L\circ{\iota},$$ as then we can identify the two properties.

For this consider the following concatenation of natural isomorphisms:
\begin{align*}
&\Ext^*_{\X\otimes\I^{op}}\left(\left(\X\otimes\I^{op}\right)_{\eta \cup 1},\iota\right)\\
&\cong \Ext^*_{\X\otimes\I^{op}}\left(\X_\eta\otimes\I^{op},\iota\right) \\
&\cong \Ext^*_{\X}\left(\X_\eta,\iota\right)\\
&\cong \Ext^*_{\X_\eta}\left(\X_\eta,L\circ \iota\right)&{\text{\cite[Corollary~5.3.2]{Rizzardo2019}}}\\
&\cong L\circ \iota.
\end{align*}

So we have that if $\widehat{\iota}$ is a colift, it immediately defines by the above isomorphisms an $\Ainfty$-lift of $L\circ \iota$ and vice versa.
\end{proof}
\end{lemma}

\begin{proposition}\label{Proposition L does not have a lift}
Let $\iota: \Inj \hookrightarrow \D^b\left(\Qcoh X\right)$ be the inclusion of injectives in the canonical heart. Then the functor $L \circ \iota: \Inj \to \D^b\left(\X_\eta\right)$ does not admit an $\Ainfty$-lift.
\begin{proof}
By Lemma~\ref{Theorem Colift and Ainfinity lifts coincide} we have that a lift of $L\circ \iota$ would be a colift of the $\Inj$-$\X$ bimodule $\iota$. This colift is obstructed by the dual characteristic morphism $c^*_{\iota,\Inj}\left(\eta\right)$. As $\M$ is invertible, we have by Lemma~\ref{Lemma comparison characteristic morphism and dual characteristic morphism} that the dual characteristic morphism vanishes if and only if the characteristic morphism vanishes. By Corollary~\ref{Corollary Inj-equivariant characteristic morphism of a concrete class is non-zero} we have that $c_{\iota,\Inj}\left(\eta\right)\neq 0.$

In particular such a colift cannot exist, and hence there also cannot exist an $\Ainfty$-lift.
\end{proof}
\end{proposition}

\section{Structure of the functor $\widetilde{f}$}\label{Section widetilde f}

It will be key for our conclusion that we carefully pick the space we consider for computing our obstructions, respectively understand the structure of the functor $\widetilde{f}$. In particular we will from now on fix the inclusion of a smooth divisor $f:X\hookrightarrow Y$, a pulled back linebundle $M=f^*\widehat{M}$ and $0\neq \eta \in \HH^m\left(X,M\right)$ with $m\geq \dim X+3$ such that $f_*\eta$. 

We start of by computing the images of sheaves $\F\in \Qcoh\left(X\right)$ under the functor $\psi_\eta$:

\begin{lemma}\label{Lemma computation image}
Let $F\in \Qcoh \left(X\right)$. Then $$\psi_\eta F\cong f_*F\oplus\left( f_*F\otimes f_*M\left[2-m\right] \right).$$
\begin{proof}
Consider the distinguished triangle in $\D^b\left(\X_\eta \right)$ induced by $L$, see \cite[Lemma~11.3]{Rizzardo2019}
$$\mathrm{can}_*\F \to L \F \to \mathrm{can}_*\Sigma^{2-m}\left(\F\otimes \M\right) \to \Sigma \mathrm{can}_*\F,$$
where we interpret $\F$ as an object in $\D^b\left(\X_\eta\right)$ via the functor $\mathrm{can}_*$.
Applying the exact functor $\widetilde{f}_*$ from Lemma~\ref{Lemma construction widetilde f} we get the following triangle in $\D^b\left(\Y\right)$
$$\left(\widetilde{f}\circ \mathrm{can}\right)_*\F \to \psi_\eta \F \to \left(\widetilde{f}\circ \mathrm{can}\right)_*\Sigma^{2-m}\left(\F\otimes \M\right) \to \left(\widetilde{f}\circ \mathrm{\can}\right)_*\Sigma \F.$$
As $\mathrm{can}\circ \widetilde{f}\cong f$ the triangle turns into
$$f_*\F \to \psi_\eta \F \to f_*\Sigma^{2-m}\left(\F\otimes \M\right) \to f_* \Sigma \F.$$
And since $f_*$ is exact we have that it is in particular graded. So we have the triangle
$$f_*\F \to \psi_\eta \F \to \Sigma^{2-m}f_*\left(\F\otimes \M\right) \to \Sigma f_*\F.$$
Now, as $\F$ and $\M$ are sheaves we have 
\begin{align*}
&\Ext^1_{\Y}\left(\Sigma^{2-m}f_*\F\otimes \M , f_*\F \right)\\
&\cong \Ext^{m-1}_{\Y}\left(f_*\F\otimes \M, f_*\F \right)\\
&\cong \Ext^{m-1}_{Y}\left(f_*\F\otimes \M , f_*\F \right)&\text{Lemma~\ref{Lemma fully faithful embedding into X}}\\
&=0 &\mathllap{\text{\cite[Lemma~11.1]{Rizzardo2019}},}
\end{align*}
where we used in the last line that $m\geq \dim X + 3$ gives $m-1 \geq \dim Y=\dim X + 1$.
So the triangle
$$f_*\F \to \psi_\eta \F \to f_*\Sigma^{2-m}\F\otimes \M \to \Sigma \F$$
splits and we have
$$\psi_\eta\left(\F\right)\cong f_*\F \oplus f_*\left(\Sigma^{2-m}M\otimes \F\right),$$
as claimed.
\end{proof}
\end{lemma}

\begin{corollary}\label{Corollary Ext concentrated in 6 degrees}
We have that the space $\Ext^*\left(\psi_\eta I, \psi_\eta J \right)$ is concentrated in degrees $m-2,m-3,1,0,3-m$ and $2-m$ for $I,J \in \Inj$.
\begin{proof}
By Lemma~\ref{Lemma computation image} we have $\psi_\eta I \cong f_*I\oplus f_*\left( I\otimes M\left[2-m\right] \right)$ which by projection formula and $M\cong f^*\widehat{M}$ turns into $f_*I\oplus  f_*\left(I\otimes \widehat{M}\right)\left[2-m\right]$. In particular we can compute
\begin{align*}
\Ext_Y^*\left(\psi_\eta I, \psi_\eta J \right)&\cong \Ext_Y^*\left(f_*I \oplus \left( f_*I \otimes \widehat{M}\right) \left[2-m\right] , f_*J\oplus \left(f_* J\otimes \widehat{M}\right)\left[2-m\right] \right)   \\
&\cong \begin{pmatrix}
\Ext_Y^*\left(f_*I,f_*J\right)& \Ext_Y^*\left(f_*I,f_* \left(J\right)\otimes \widehat{M}[2-m]\right)\\
\Ext_Y^*\left(f_* \left(I\right)\otimes \widehat{M}[2-m],f_* J\right)&\Ext_Y^*\left(f_* I,f_* J \right)
\end{pmatrix}
 \end{align*}
Where we applied in the bottom right corner the fact that $\widehat{M}$ is invertible while also dropping the shifts.
Now it suffices to observe that we have
\begin{align*}
\Ext_Y^*\left(f_*F,f_*G\right) &\cong  \Ext^{*}_X\left(f^*f_*F,G\right) \\
&\cong \Ext^{*}_X\left(I\otimes C_f[-1]\oplus F,G\right) & \mathllap{f^*f_* F \cong \left(F \oplus F \otimes C_f[1]\right)}\\
&\cong \Ext^{*}_X\left(F,G\right)\oplus\Ext^*_X\left(F\otimes C_f,G\right) &\quad \quad \Ext^*\text{ is additive} 
\end{align*}
Where we use $C_f$ for the conormal bundle of the inclusion $f:X\hookrightarrow Y$. Now it suffices to observe that as $f$ is a divisor $C_f$ is invertible too and in our cases $F$ and $G$ are either injective, or an injective tensored by an invertible module $M$.
\end{proof}
\end{corollary}

In order to have control over the obstructions explicitly, we will study the functor $\widetilde{f}:\Y \to \X_\eta$ from Lemma~\ref{Lemma construction widetilde f} explicitly. For this consider the commutative diagram of $\Ainfty$-categories and functors giving rise to a candidate.

\begin{equation}\label{Diagram widetilde f explicitly}
\tikz[heighttwo,xscale=3,yscale=3,baseline]{
\node (X) at (0,0) {$\X$};
\node (Xeta) at (0,1) {$\X_\eta$};
\node (Yfeta) at (1,1) {$\Y_{f_*\eta}$};
\node (Y0) at (2,1) {$\Y_0$};
\node (YM) at (3,1) {$\Y \oplus f_*\M[2-m]$};
\node (Y) at (2,0) {$\Y$};
\draw[<-]
(Xeta) edge node[above] {$f$} (Yfeta)
(X) edge node[above] {$f$} (Y)
(Yfeta) edge node[above] {$h$} (Y0);
\draw[->]
(Xeta) edge node[right] {$\text{can}$} (X)
(Yfeta) edge node[above right] {$\text{can}$} (Y)
(Y0) edge node[right] {$\text{can}$} (Y)
(YM) edge node[below right] {$\text{can}=\left(\id,0\right)$} (Y)
;
\draw[->, dashed]
(Y) edge node[above right] {$\widetilde{f}$} (Xeta)
;
\draw[<-]
(Y0) edge node[above] {$=$} (YM);
}
\end{equation}
The bottom arrow is the functor given in Lemma~\ref{Lemma associated f between k-linear categories} and we denote by $\text{can}$ the canonical functors arising from the infinitesimal deformation, see Definition~\ref{Definition infinitesimal deformation along Hochschild cohomology class}. 

The other functors are given by:
\begin{itemize}
\item The functor  $f: \Y_{f_*\eta}\to \X_\eta$ is the functor given by the identity on objects and acting on morphisms by 
\begin{align*}
\begin{pmatrix}
f&0\\
0&f
\end{pmatrix}:\Y\left(i,j\right)\oplus f_*\M\left(i,j\right)[2-m] \to \X\left(i,j\right)\oplus \M\left(i,j\right)[2-m].
\end{align*}
This functor is by construction an $\Ainfty$-functor as by definition $f_*\eta$ is given by the composition of bimodule morphisms
$$f_*\eta: \Y \xrightarrow{f} \X \xrightarrow{\eta}\M[2-m] .$$ 

\item The functor $h: \Y_0 \to \Y_{f_*\eta}$ is given by a nullhomotopy of $$0=f_*\eta \in \HH^*\left(\Y,f_*M\right)$$ from Proposition~\ref{Proposition construction X eta independent of choice of cocycle} and is in particular an $\Ainfty$-equivalence.
\item The functor $\Y \oplus \M[2-m]\to\Y_0 $ is the definition of the trivial infinitesimal deformation.
\end{itemize}

The diagram \eqref{Diagram widetilde f explicitly} now shows that the main choice involved in the definition of $\widetilde{f}$ is the choice of a split $s: \Y \to \Y\oplus f_*\M[2-m]$. For which we take the one given by 
$$s=\begin{pmatrix}
\id\\0
\end{pmatrix}: \Y \to \Y\oplus f_*\M[2-m].$$

In particular we can use these constructions to study the functors involved in the construction of $\widetilde{f}$ explicitly and conclude that we can restrict to a suitable space where we can control the $\Ainfty$-bimodule structure rigorously.

\begin{definition}\label{definition widetilde f by composition}
We consider for the remainder of this paper the concrete choice of the functor $\widetilde{f}: \Y \to \X_\eta$ given by the composition
$$\widetilde{f}: \Y \xrightarrow{s}\Y_0\xrightarrow{h}\Y_{f_*\eta}\xrightarrow{f}\X_\eta.$$
\end{definition}

Now we continue our study by restricting to the category of injectives in $\D^b\left(\Qcoh X\right)$, $$\iota: \Inj \hookrightarrow \D^b\left(\Qcoh X\right).$$

For this consider the following composition of all intermediate steps of the construction of $\psi_\eta$. We also include a step for the essential image of $\Inj$ in $\D^b\left(\X_\eta\right)$ respectively in $\D^b\left(\Y\right)$, which we denote by $\Inj_\eta$, respectively $f_*\Inj_{f_*\eta}$,

\begin{equation}\label{Diagram widetilde f explicitly for Injectives}
\tikz[heighttwo,xscale=2.5,yscale=2.5,baseline]{
\node (Inj) at (0,0) {$\Inj$};
\node (DX) at (1,0) {$\D^b\left(\X\right)$};
\node (Xeta) at (0,1) {$\Inj_\eta$};
\node (fXfeta) at (1,1) {$f_*\Inj_{f_*\eta}$};
\node (Yfeta) at (2,1) {$\D^b\left(\Y \right)_{f_*\eta}$};
\node (YM) at (3,1) {$\D^b\left(\Y_0\right) $};
\node (Y) at (3,0) {$\D^b\left(\Y\right).$};
\draw[->]
(Xeta) edge node[above] {$f_*$} (fXfeta)
(DX) edge node[above] {$\psi_\eta$} (Y)
(Yfeta) edge node[above] {$h_*$} (YM)
(YM) edge node[right] {$s_*$} (Y)
(Inj) edge node[left] {$L\circ \iota$}(Xeta);
;
\draw[right hook ->]
(Inj) edge node[above] {$\iota$}(DX)
(fXfeta) edge node[above]{$j$}(Yfeta);
}
\end{equation}

\begin{remark}
The diagram $\eqref{Diagram widetilde f explicitly for Injectives}$ allows us to point out our proof strategy more clearly. As we want to obstruct the composition $\psi_\eta \circ \iota$ from having an $\Ainfty$-lift, we will prove that the composition of $h_*$ and $s_*$ does not change the obstruction and so we may restrict ourselves to the category $f_*\Inj_{f_*\eta}$ for the computation of obstructions. We then use the explicit structure of the category $f_*\Inj_{f_*\eta}$ to compute that these obstructions cannot vanish independent of intermediate lifts.

Conceptually one can think about this argument as the bimodule $\psi_\eta$ not being $\Ainfty$ over $\Inj$ while being $\Ainfty$ over $\Y$. In particular it behaves like a non-$\Ainfty$ Fourier-Mukai kernel.
\end{remark}

The following Proposition~\ref{Proposition structure essential image of L iota} and Proposition~\ref{Proposition structure essential image of f* L iota} show that the notations $\Inj_\eta$ and $f_*\Inj_{f_*\eta}$ for the essential image of $f_*\circ L\circ \iota$ are warranted.

\begin{proposition}\label{Proposition structure essential image of L iota}
The essential image of $L\circ \iota$ has objects $I\in \Inj$, morphism spaces given by 
$$\textrm{essim}\left(L\circ \iota\right)\left(I,J\right) \cong \Ext_X^{0}\left(I,J\right)\oplus \Ext_X^{0}\left(I,J\otimes \M\right)[2-m]$$
and $\m_m=\eta$.
\begin{proof}
By construction of $L$, see \cite[(5.7)]{Rizzardo2019}, the morphism spaces in $\Inj_\eta$ are of the shape 
\begin{align*}
\mathrm{essim}\left(L\circ\iota\right)\left(I,J\right)&\cong \Ext_X^{*}\left(I,J\right)\oplus \Ext_X^{2-m+*}\left(I,J\otimes \M\right) \\&\cong \Ext_X^{0}\left(I,J\right)\oplus \Ext_Y^{0}\left(I,J\otimes \M\right)[2-m].
\end{align*}
As $L$ is defined by identifying the inclusion of $\Inj$ into $\D^b\left(\X\right)$ with the inclusion into $\D^b\left(\X_\eta\right)$ \cite[Lemma~6.3.1]{Rizzardo2019} we have that there is a non-trivial $\m_m=\eta$, via the identification from \cite[Theorem~8.17]{Lowen2005} $\HH^*\left(\Qcoh,\M\right)\cong \HH^*\left(\Inj,\M\right)$. 
\end{proof}
\end{proposition}

\begin{proposition}\label{Proposition structure essential image of f* L iota}
The essential image of $f_*\circ L\circ \iota$ has objects $I\in \Inj$ and morphism spaces given by 
$$\textrm{essim}\left(f_*\circ L\circ \iota\right)\left(I,J\right) \cong \Ext_Y^{*}\left(f_*I,f_*J\right)\oplus \Ext_Y^{2-m+*}\left(f_*I,f_*\left(J\otimes \M\right)\right).$$
\begin{proof}
As $f_*$ is an $\Ainfty$-functor it suffices to push forward the structure from Proposition~\ref{Proposition structure essential image of L iota}.
\end{proof}
\end{proposition}

\begin{proposition}\label{proposition s* is inclusion of upper triangular matrices}
The functor $s_*:\D^b\left(Y_0\right)\to \D^b\left(\Y\right)$ acts on objects of the shape $h_*\circ f_* \circ L \left(I\right)$ for $I\in Inj$ by
\begin{align*}
s_*:\begin{pmatrix}
\Ext_Y^*\left(f_*I,f_*J\right)\\ \Ext_Y^{*+m-2}\left(f_*I,f_*\left(J\otimes M\right)\right)\end{pmatrix} &\to \Ext_Y^*\left(f_*I\oplus \Sigma^{m-2}f_*\left(I\otimes M \right),f_*J\oplus \Sigma^{m-2} f_*\left(J\otimes M\right) \right)\\
\begin{pmatrix}\alpha\\ \beta\end{pmatrix} &\mapsto \begin{pmatrix}
\alpha & \beta \\ 0&\alpha \end{pmatrix}
\end{align*}
\begin{proof}
By the above discussions and Proposition~\ref{Proposition structure essential image of f* L iota} we have 
$$\mathrm{essim}\left(f_*\circ L \circ \iota\right)\left(I,J\right)\cong \Ext_Y^{*}\left(f_*I,f_*J\right)\oplus \Ext_Y^{2-m+*}\left(f_*I,f_*\left(J\otimes \M\right)\right).$$
As $h$ is a homotopy equivalence the same holds for $\mathrm{essim}\left(h_*\circ f_*\circ L \circ \iota\right)$.
By Lemma~\ref{Lemma computation image} we have $\psi_\eta \left(I\right)\cong f_*I\oplus f_*\left(M\otimes I\right)[2-m]$. As we made the choice for $s$ being the trivial inclusion $\Y\hookrightarrow \Y_0\cong \Y \oplus f_*\M[2-m]$, $s_*$ induces the inclusions of upper triangular matrices with identical diagonal as claimed.
\end{proof}
\end{proposition}

\begin{corollary}\label{Lemma the functor s* is an equivalence on negative ext}
The functor $s_*:\D^b\left(\Y_0\right)\to \D^b\left(\Y \right)$ is an isomorphism on negative degree morphisms in the essential image of $\Inj$.
\begin{proof}
By Proposition~\ref{proposition s* is inclusion of upper triangular matrices} we have that in the matrix expression of $$\begin{pmatrix}
\alpha & \beta \\ \gamma&\delta \end{pmatrix} \in \Ext_Y^*\left(f_*I\oplus \Sigma^{m-2}f_*\left(I\otimes M \right),f_*J\oplus \Sigma^{m-2} f_*\left(J\otimes M\right) \right)$$ 
the morphisms $\alpha,\delta$ and $\gamma$ must be in non-negative degrees as the direct summands $\Ext^*_Y\left(f_*I,f_*J\right)$,$\Ext_Y^*\left(f_*I\otimes M, f_*J\otimes M\right)$ and $\Ext_Y^*\left(\Sigma^{m-2}f_*I\otimes M,f_*J\right)$ are concentrated in non-negative degrees.

 In particular it suffices to verify the claim for $\beta$. There we have that the map 
\begin{align*}
s_*:\Ext_Y^{*+m-2}\left(f_*I,f_*\left(J\otimes M\right)\right) &\to \Ext_Y^*\left(f_*I\oplus \Sigma^{m-2}f_*\left(I\otimes M \right),f_*I\oplus \Sigma^{m-2} f_*\left(I\otimes M\right) \right)\\
\beta &\mapsto \begin{pmatrix}
0 & \beta \\ 0&0 \end{pmatrix}
\end{align*}
is an isomorphism given by the inclusion of the direct summand
$$\Ext_Y^{*}\left(f_*I,\Sigma^{m-2}f_*\left(J\otimes M\right)\right) \hookrightarrow \Ext_Y^*\left(f_*I\oplus \Sigma^{m-2}f_*\left(I\otimes M \right),f_*I\oplus \Sigma^{m-2} f_*\left(I\otimes M\right) \right).$$
After restricting to the negative morphisms this remains an isomorphism, which finishes the claim.
\end{proof}
\end{corollary}

\section{Computing obstructions explicitly for Divisors}\label{Section Computing obstructions explicitly}

The observations from the previous sections let us compute:

\begin{proposition}\label{Proposition obstruction of L does not vanish}
The obstruction $$o_m\left(L\circ \iota \right)=[m_m,L]\in \HH^*\left(\Inj, \Inj_\eta\right)_{2-m}$$
does not vanish.
\begin{proof}
We know by Proposition~\ref{Proposition L does not have a lift} that $L\circ \iota$ does not have an $\Ainfty$-lift. In particular we have by Lemma~\ref{Lemma A-infinity obstructions against lift of functors} that there must be a non-vanishing 
$$o_i\left(L\circ \iota \right) \in \HH^i\left(\Inj, \Inj_\eta\right)_{i}.$$ 
However, as $L\circ \iota$ is a functor it is $\A_2$ and the non-vanishing obstruction must have $i\geq 3$. 
As we have by Proposition~\ref{Proposition structure essential image of L iota} that $\Inj_\eta$ is concentrated in degrees $0$ and $2-m$ we have that the only non-trivial obstruction can be
$$o_m\left(L\circ \iota \right)\in \HH^*\left(\Inj, \Inj_\eta\right)_{2-m}$$
and in particular $o_m\left(L\circ \iota \right)\neq 0$.
\end{proof}
\end{proposition}

\begin{lemma}\label{Lemma comparison of obstructions}
The obstructions of $\A_i$ functors from a $\kk$-linear category with essential image $\mathrm{essim}\left(\psi_\eta\right)$ can be computed on $f_*\Inj_{f_*\eta}$, c.f. Lemma~\ref{Lemma A-infinity obstructions against lift of functors}.

More precisely let $$\varphi=\left(\varphi_1,...,\varphi_i\right): \I\to f_*\Inj_{f_*\eta}$$ be an $\A_i$ functor from a $\kk$-linear category $\I$ then we have an isomorphism
\begin{align*}
\left(s_*\circ h_*\circ j\right)_* : \HH^*\left(\I,f_*\Inj_{f_*\eta}\right)_{2-m} &\xrightarrow{\sim}\HH^*\left(\I,\D^b\left(\Qcoh Y\right)\right)\\
o_i\left(f_*\circ \varphi \right) &\mapsto o_i\left(s_*\circ h_* \circ j \circ \varphi\right).
\end{align*}
Furthermore an $\A_{i+1}$-lift $\varphi_{i+1}$ of $s_*\circ h_* \circ j \circ \varphi  $ gives rise to an $\A_{i+1}$-functor $\widehat{\varphi_{i+1}}$ such that $$\varphi_{i+1} =s_*\circ h_* \circ j \circ \widehat{\varphi}_{i+1}.$$

\begin{proof}
As $h$ is a homotopy equivalence and $j$ is the inclusion of a full subcategory we have by naturality of obstructions 
\begin{align*}
\left(h\circ j\right)_*: \HH^m\left(\I,f_*\Inj_{f_*\eta}\right) &\xrightarrow{\sim} \HH^m\left(\I, \D^b\left(\Y_0\right)\right)\\
o_m\left(\varphi \right)&\mapsto o_m\left(h_*\circ j \circ \varphi \right).
\end{align*}
For the compatibility with lifts observe that since $j$ and $h$ are homotopy equivalences on the image of the injectives they admit quasi-inverses they are immediately compatible with intermediate lifts.

As by Proposition~\ref{proposition s* is inclusion of upper triangular matrices} the functor $s_*$ is the inclusion of the upper triangular matrices and by Proposition~\ref{Lemma the functor s* is an equivalence on negative ext} an equivalence on the negative $\Ext$-groups. We get, since $\I$ is a $\kk$-linear category that $s_*$ is an equivalence,
\begin{align*}
s_*:\HH^m\left(\I, \D^b\left(\Y_0\right)\right)_{< 0}&\xrightarrow{\sim} \HH^m\left(\I,\D^b\left(\Qcoh Y\right)\right)_{< 0}\\
o_m\left(h_*\circ j \circ \varphi \right)&\mapsto o_m\left( s_*\circ h_*\circ j \circ \varphi \right).
\end{align*}
Since any intermediate lift would land in the $\Ainfty$-submodule of morphisms in negative degrees and $s_*$ is an isomorphism on these we have that also this comparison morphism is compatible with intermediate lifts. 

So altogether we have the desired isomorphism and compatibility with lifts.
\end{proof}
\end{lemma}

The above Lemma~\ref{Lemma comparison of obstructions} has the following immediate consequence

\begin{corollary}\label{Corollary comparison Ainfty lifts}
The functor $\psi_\eta \circ \iota:\Inj \to \D^b\left(\Qcoh Y\right)$ has an $\Ainfty$-lift if and only if $f_*\circ L \circ \iota: \Inj \to f_*\Inj_{f_*\eta}$ has an $\Ainfty$-lift.
\end{corollary}

Now we may restrict for the computation of the obstructions $$o_i\left(\psi_\eta\circ \iota \right)\in \HH^i\left(\Inj,\H^{2-i}\D^b_\dg\left(\Y\right)\right)$$ to the image of $f_*\circ L$ in negative degrees. In particular, as we consider minimal $\Ainfty$-structures, the $\Inj$-bimodule structure on 
\begin{align*}
\H^{<0}\mathrm{im}\left(f_*\circ L \circ \iota \right)&\cong \Ext^{*-2+m}\left(f_*I,f_*\left(I\otimes M\right)\right)&\text{Corollary~\ref{Corollary Ext concentrated in 6 degrees}}
\end{align*}
 is given by the $\kk$-linear functor
\begin{align*}
\H^{<0}\mathrm{im}\left(f_*\circ L \circ \iota \right):\Inj \otimes \Inj^{\op}&\to \D^b_\infty \left(\kk\right)\\
\left(I,J\right)&\mapsto \Ext^{*-2+m}_Y\left(f_*I,f_*J\otimes M\right)\\
\left(\alpha\otimes \beta\right)&\mapsto f_*\circ L\left(\alpha\right)\otimes \left(f_*\circ L \left(\beta\right)\otimes f_*M\right).
\end{align*}

\begin{lemma}\label{Lemma f* is an equivalence in degree 2-m}
The functor $f_*: \Inj_\eta \to f_*\Inj_{f_*\eta}$ is an equivalence on $\H^{2-m}\Inj_\eta$.
$$f_*: \H^{2-m}\Inj_\eta \xrightarrow{\sim} \H^{2-m}f_*\Inj_{f_*\eta}$$
\begin{proof}
We know that $f_*$ is by construction literally the functor $f_*$ associated to the closed immersion of a divisor $f:X\hookrightarrow Y$ on the two components of 
$$\Inj_\eta \left(I,J\right)\cong \Ext^0_X\left(I,J\right)\oplus \Ext^0_X\left(I,J\otimes M\right)[2-m].$$
Now, as $f_*$ is induced by a direct image along a closed immersion it is an isomorphism on $\Ext^0$-groups. In particular  
\begin{align*}
f_*: \H^{0}\Inj_\eta &\xrightarrow{\sim} \H^{0}f_*\Inj_{f_*\eta}\\
f_*: \H^{2-m}\Inj_\eta &\xrightarrow{\sim} \H^{2-m}f_*\Inj_{f_*\eta}
\end{align*}
as claimed.
\end{proof}
\end{lemma}

\begin{proposition}\label{Proposition essential image of L is an Inj A infinity bimodule}
The functor $\H^{<0}\mathrm{im}\left(f_*\circ L \circ \iota \right)$ is naturally an $\Inj$-$\Ainfty$-bimodule with action given via the functor $f_*\circ L$.
\begin{proof}
We have that $\Inj$ is $\kk$-linear and $\H^{<0}\mathrm{im}\left(f_* \circ L \circ \iota \right)$ is concentrated in degrees $2-m$ and $3-m$. So the only possible higher multiplications of $\Inj$ on $\H^{<0}\mathrm{im}\left(f_* \circ L \circ \iota \right)$ are $\m_2$ and $\m_3$. By construction we have that $f_*\circ L$ is a functor, so it is compatible with $\m_2$ and in particular associative, i.e.
$$\left[\m^{f_*\circ L}_2,m^{f_*\circ L}_2\right]=0.$$ 

Now for $\m_3$ we need to verify 
$$\left[\m^{f_*\circ L}_3,\m^{f_*\circ L}_2\right]=0.$$
This is by minimality the $\Ainfty$-equation for the whole category $$\H^*\mathrm{im}\left(f_* \circ L \circ \iota \right)\cong f_*\Inj_{f_*\eta}.$$
In particular 
$$\left[\m^{f_*\Inj_{f_*\eta}}_3,\m^{f_*\Inj_{f_*\eta}}_2\right]=\left[\m^{f_*\circ L}_3,\m^{f_*\circ L}_2\right]=0$$ holds.

\end{proof}
\end{proposition}

Before we continue to the computation of $\m_3$ on $\H^{<0}\mathrm{im}\left(f_* \circ L \circ \iota \right)$ as an $\Inj$-bimodule we recall the following fact.

\begin{lemma}\label{Lemma A-infty adjunction is morphism of bimodules}
Let $f^*\dashv f_*:\D^b_\infty \left(\Qcoh X\right)\to \D_\infty^b \left(\Qcoh Y\right)$ be an adjuction of $\Ainfty$-functors and let $N$ be an $\Ainfty$-$\I$-bimodule of the shape
\begin{align*}
N: \I \otimes \I^\op &\to \D^b_\infty\left(\Qcoh X\right)\\
\left(i,j\right)&\mapsto \Ext^*_X\left(Ni,Nj\right)\\
\alpha \otimes \beta &\mapsto  N_1 \left(\alpha\right)\otimes N_2\left(\beta\right)
\end{align*} 
for $\I$ a $\kk$-linear category. Then we have that the adjunction quasi-isomorphism 
$$\psi:\Ext^*_\B\left(f_*X,f_*Y\right)\xrightarrow{\sim}\Ext^*_\A\left(f^*f_*X,Y\right)$$
induces a quasi-isomorphism of $\Ainfty$-$\I$-bimodules.

$$
\begin{pmatrix}
f_* N: \I \otimes \I^\op \to \D^b\left(\Qcoh Y\right)\\
\left(i,j\right)\mapsto \Ext^*_X\left(f_*Ni,f_*Nj\right)\\
\alpha \otimes \beta \mapsto  f_* N_1 \left(\alpha\right)\otimes f_* N_2\left(\beta\right)
\end{pmatrix}
\cong \begin{pmatrix}
 \I \otimes \I^\op \to \D^b\left(\Qcoh X\right)\\
\left(i,j\right)\mapsto \Ext^*_X\left(f^*f_* Ni,Nj\right)\\
\alpha \otimes \beta \mapsto  f^*f_* N_1 \left(\alpha\right)\otimes N_2\left(\beta\right)
\end{pmatrix}.
$$
\begin{proof}
Consider a $\dg$-replacement of the adjunction $$f^*\dashv f_*:\D^b_\dg\left(\Qcoh X\right)\to \D^b_\dg\left(\Qcoh Y\right).$$
Then it suffices to show that the isomorphism
$$\D^b_\dg\left(\Qcoh Y \right)\left(f_*N i,f_* N j\right)\xrightarrow{\sim}\D^b_\dg\left(\Qcoh X\right)\left(f^*f_*Ni,Nj\right)$$
is a quasi-isomorphism of $\dg$-$\I$-bimodules.
By definition of $\dg$-adjunction it is already a quasi-isomorphism of chain complexes of vectorspaces. In particular it suffices to prove that it is a morphism of $\I$-bimodules.

For this consider first the operation from the left on $f_*N$. Denote it by 
$$f_*N_1\left(\alpha\right):f_*Ni'\to f_*Ni ,$$
then we have the commutative diagram

$$\tikz[heighttwo,xscale=6,yscale=3,baseline]{
\node (DX) at (0,0) {$\D_\dg\left(\Qcoh Y\right)\left(f_*N i',f_*Nj\right)$};
\node (DQX) at (0,1) {$\D_\dg\left(\Qcoh Y\right)\left(f_*N i,f_*Nj\right)$};
\node (DQY) at (1,1) {$\D_\dg\left(\Qcoh X\right)\left(f^*f_*Ni,Nj\right)$};
\node (DY) at (1,0) {$\D_\dg\left(\Qcoh X\right)\left(f^*f_*Ni',Nj\right).$};

\draw[->]
(DX) edge node[above] {$\sim$} (DY)
(DQX) edge node[left] {$\_\circ f_*N_1\left(\alpha\right)$} (DX)
(DQY) edge node[right] {$\_\circ f^*f_*N_1\left(\alpha\right)$} (DY)
(DQX) edge node[above] {$\sim$} (DQY)
;
}
$$

Now for the action on the left given by $f_*N_2\left(\beta\right):f_*Nj\to f_*Nj'$ we get analogously

$$\tikz[heighttwo,xscale=6,yscale=3,baseline]{
\node (DX) at (0,0) {$\D_\dg\left(\Qcoh Y\right)\left(f_*Ni,f_*Nj'\right)$};
\node (DQX) at (0,1) {$\D_\dg\left(\Qcoh Y\right)\left(f_*Ni,f_*Nj\right)$};
\node (DQY) at (1,1) {$\D_\dg\left(\Qcoh X\right)\left(f^*f_*Ni,Nj\right)$};
\node (DY) at (1,0) {$\D_\dg\left(\Qcoh X\right)\left(f^*f_*Ni,Nj'\right).$};

\draw[->]
(DX) edge node[above] {$\sim$} (DY)
(DQX) edge node[left] {$\_\circ f_* N_2\left(\beta\right)$} (DX)
(DQY) edge node[right] {$\_\circ N_2\left(\beta\right)$} (DY)
(DQX) edge node[above] {$\sim$} (DQY)
;
}
$$
In particular 
$$\D^b_\dg\left(\Qcoh\left(Y\right)\right)\left(f_*Ni,f_*Nj\right)\xrightarrow{\sim}\D^b_\dg\left(\Qcoh X\right)\left(f^*f_*Ni,Nj\right)$$
is compatible with both actions and induces a quasi-isomorphism of $\Ainfty$-$\I$-bimodules.
\end{proof}
\end{lemma}

\begin{remark}
In the above proof it is essential that we only consider the bimodule structure given by the functor $f_*$. Otherwise we would not be able to compare the right action, as not every functor in the essential image arises via $f_*$.
\end{remark}

The following is a rewritten proof to our setup of the well known fact that for a divisor $X\hookrightarrow Y$ and a sheaf $\F$ the  object $f_* \F \in \D^b_\dg\left(Y\right)$ usually is not formal over $Y$, however, it is formal over $X$. In particular it is essential that we restrict in this proof to the bimodule structure over $\Inj$.

\begin{proposition}\label{Proposition Himpsi is formal on injectives}
The  $\Inj$-$\Ainfty$-bimodule $\H^{<0}\mathrm{im}\left(f_* \circ L \circ \iota \right)$ has  $\m_3=0$.
\begin{proof}

By adjunction we have a natural $\Ainfty$-bimodule-isomorphism written down objectwise $$\Ext_Y^{*-2+m}\left(f_*I,f_* J \otimes M\right) \cong \Ext^{*-2+m}_X\left(f^*f_*I,J\otimes M\right).$$
This isomorphism is by Lemma~\ref{Lemma A-infty adjunction is morphism of bimodules} compatible with the bimodule structure induced by $\X$ on the $\X$-bimodule $\H^{<0}\mathrm{im}\left(f_*\circ L \circ \iota \right)$ and so in particular compatible with the $\Inj$-bimodule structure induced by $f_*\circ L \circ \iota$ from Proposition~\ref{Proposition essential image of L is an Inj A infinity bimodule}.

So we have the following chain of isomorphisms of $\Inj$-$\Ainfty$-bimodules, written down objectwise:
\begin{align*}
&\H^{<0}\mathrm{im}\left(f_* \circ L \circ \iota \right)\left(I,J\right) \\
&\cong \Ext^{*-2+m}_X\left(f^*f_*I,J\otimes M\right) & f^* \vdash f_*\\
&\cong \Ext^{*-2+m}_X\left(f^*f_*\OO_X \otimes_{\OO_Y}I,J\otimes M\right) \\
&\cong \Ext^{*-2+m}_X\left(I\otimes C_f [1]\oplus I,J\otimes M\right) & \mathllap{f^*f_*\OO_X \cong \otimes \left(\OO_X \oplus C_f[1]\right)}\\
&\cong \Ext^{*-2+m}_X\left(I,J\otimes M\right)\oplus\Ext^{*-3+m}_X\left(I\otimes C_f,J\otimes M\right) &\Ext^*\text{ is additive}
\end{align*}
where all isomorphisms hold on the $\Ainfty$-level. The first one is an isomorphisms of $\Ainfty$-$\Inj$-bimodules by Lemma~\ref{Lemma A-infty adjunction is morphism of bimodules}, the second is the defintition of the pullback, the third one is $\Ainfty$/$\dg$ by \cite[0.7]{Arinkin2012} and the last one is $\Ainfty$ as it is induced by the additivity of $\Ext$.
Altogether this sequence of natural isomorphisms operates on actions via:
\begin{align*}
f_* \circ L \left(\alpha\right)&\mapsto 
\begin{pmatrix}
 L \left(\alpha\right) \\
  L \left(\alpha\right) \otimes C_f\end{pmatrix}\\
 f_* \circ L \left(\beta\right)&\mapsto   \begin{pmatrix} L \left(\alpha\right) \\
  L \left(\alpha\right) \end{pmatrix} .
\end{align*}

As $\Inj$ is $\kk$-linear and both direct summands are concentrated in one degree we have that the $\m_3$ on $\Ext^*_X\left(I,J\otimes M\right)$ and $\Ext^*_X\left(I\left(-X\right),J\otimes M\right)$ vanishes. In particular we found a minimal $\Ainfty$-structure with vanishing $\m_{3}$ of the $\Inj$-$\Ainfty$-bimodule
$$\H^{<0}\mathrm{im}\left(f_* \circ L \circ \iota \right)\cong \H^{3-m}\mathrm{im}\left(f_* \circ L \circ \iota \right)\oplus \H^{2-m}\mathrm{im}\left(f_* \circ L \circ \iota \right)$$
as claimed. 
\end{proof}
\end{proposition}

\begin{remark}
Even though we just proved that $\m_3=0$ on the $\Ainfty$-bimodule $\H^{<0}\mathrm{im}f_*\circ L\circ \iota$, this does not apply to $\m_m=f_* \eta$ on $\H^*\mathrm{im}f_*\circ L\circ \iota$. Even more, due to $\m_{\geq 3}$ vanishing on the $\Inj$-bimodule $\H^{\le 0}\mathrm{im}f_*\circ L\circ \iota$, we actually have that  $\mathrm{im}f_*\circ L$ has a non-vanishing $\m_m=f_*\eta$ as an $\Inj$-bimodule.
\end{remark}

\begin{lemma}\label{Lemma obstruction is independent and non-zero intermediate}
We have independent of choices $$0\neq o_m\left(f_*\circ L \circ \iota\right)\in \HH^*\left(\Inj, \H^{<0}\mathrm{\im}f_*\circ L\circ \iota \right)$$
\begin{proof}
By naturality of obstructions we have 
\begin{align*}
o_{m-1}\left(f_*\circ L \circ \iota\right)&=f_*o_{m-1}\left(L \circ \iota\right)=0
\end{align*}
 In particular we may choose an $\A_{m-1}$ lift $$\left(f_*\circ L \circ \iota, 0,..., 0, \delta, g \right): \Inj\to \H^{\le 0} \mathrm{\im}f_*\circ L \circ \iota \cong f_*\Inj_{f_*\eta}.$$

By Proposition~\ref{Proposition Himpsi is formal on injectives} we have
 $$\m_3 = 0$$
  on the $\Inj$-bimodule $\H^{<0}\mathrm{\im}f_*\circ L\circ \iota$. Now we compute in 
$$\HH^m\left(\Inj, \H^{2-m}\mathrm{\im}\psi_\eta\circ \iota \right)=\HH^m\left(\Inj, \H^{2-m}\mathrm{\im}^{< 0}\psi_\eta\circ \iota \right)$$ that the obstruction is independent of choices
\begin{align*}
o_m\left(f_*\circ L \circ \iota, 0,..., 0, \delta, g \right)&= [\m_m,f_*\circ L \circ\iota] \pm [\m_3,\delta] \pm [\m_2,g] &\eqref{Equation defining obstructions}\\
&=[\m_m,f_*\circ L \circ\iota] \pm [\m_3,\delta]  &[\m_2,\_]=\mathrm{d}_{\mathrm{Hoch}}\\
&=[\m_m,f_*\circ L \circ\iota]  & \text{Proposition~\ref{Proposition Himpsi is formal on injectives}}\\
&=\m_m\left(f_*\circ L \circ \iota ,...,f_*\circ L \circ \iota  \right)  &\Inj \text{ has }\m_m=0\\
&=f_* \m_m\left(L\circ\iota,...,L \circ \iota \right) & f_* \text{ is }\Ainfty\\
&=f_* o_m\left(L\circ \iota \right)\\
&\neq 0.
\end{align*}
Where we used in the last line that $f_*$ induces an equivalence by Lemma~\ref{Lemma f* is an equivalence in degree 2-m} on 
$$\Ext_{\X_\eta}^{2-m}\left(L\left(I\right),L\left(J\right)\right)\cong \Ext_{\Y_{f_*\eta}}^{2-m}\left(f_*\circ L\left(I\right),f_*\circ L\left(J\right)\right).$$
In particular the obstruction does not vanish and is independent of choices.
\end{proof}
\end{lemma}

\begin{proposition}\label{Proposition obstruction is independent and non-zero}
We have independent of choices $$0\neq o_m\left(\psi_\eta \circ \iota\right).$$
\begin{proof}
By Lemma~\ref{Lemma comparison of obstructions} we have have
$$o_m\left(\psi_\eta \circ \iota\right)=o_m\left(f_*\circ L\circ \iota\right)$$
and by Lemma~\ref{Lemma obstruction is independent and non-zero intermediate} we have independent of choices 
$$o_m\left(f_*\circ L\circ \iota\right)\neq 0.$$
Which finishes the claim.
\end{proof}
\end{proposition}

\begin{theorem}\label{Theorem non-Fourier-Mukai Qcoherent}
The functor $$\psi_\eta:\D^b\left(\Qcoh X \right)\to \D^b\left(\Qcoh Y\right)$$ does not admit an $\Ainfty$/$\dg$-lift.
\begin{proof}
By Lemma~\ref{Proposition obstruction is independent and non-zero} we have that $o_m\left(\psi_\eta \circ \iota\right)$ does not vanish and is independent of choices. In particular $\psi_\eta \circ \iota$ cannot have an $\Ainfty$-lift. Since $\iota$ is dg by Lemma~\ref{Lemma inclusion of injectives is dg}, this means that $\psi_\eta$ cannot have an $\Ainfty$-lift. 
\end{proof}
\end{theorem}

If the functor $\psi_\eta$ is continuous we can restrict our argument to the category of coherent sheaves. We conjecture that $\psi_\eta$ has this property as it is by definition defined on the subcategory of injectives and extended from there, see \cite{Rizzardo2019}. However, as it is very ill-behaved we were not yet able to verify it.

\begin{theorem}\label{Theorem non-Fourier-Mukai Coherent}
Assume that $\psi_\eta$ is continuous. Then we have that the restriction $$\psi_\eta: \D^b\left(X\right)\to \D^b\left(Y\right)$$ is not Fourier-Mukai.
\begin{proof}
Assume that $\psi_\eta: \D^b\left(X\right)\to \D^b\left(Y\right)$ admits a dg-enhancement. Let $\widehat{\psi}_\eta: \D^b\left(\Qcoh X\right) \to \D^b\left(\Qcoh Y\right)$ be the continuation of $\varphi_\eta$ to $\D^b\left(\Qcoh \_ \right)$ via \cite[Proposition~5.3.5.10.]{Lurie}, respectively \cite[Theorem~7.2]{Toen2004}, we have that this is the unique continuous extension of $\psi_\eta$. In particular we have, as $\psi_\eta$ was assumed to be continuous 
$$\widehat{\psi_\eta}=\psi_\eta::\D^b\left(\Qcoh X \right)\to \D^b\left(\Qcoh Y\right).$$
But as $\psi_\eta: \D^b\left(X\right)\to \D^b\left(Y\right)$ had a $\dg$-enhancement, $\widehat{\psi_\eta}$ would have a $\dg$-enhancement, which is contradiction to Theorem~\ref{Theorem non-Fourier-Mukai Qcoherent}
\end{proof}
\end{theorem} 

This in particular means that the following corollaries could be restricted to the case of coherent sheaves if $\psi_\eta$ is continuous.

By \cite[Theorem~4.14]{Kueng2022} we can compute for a smooth hypersurface the kernels of $\HH^*\left(X,\OO\left(p\right)\right)$ explicitly via twisted Hodge-diamonds, in particular all these kernels in  $\HH^{\geq \dim X + 3}\left(X,\OO\left(p\right)\right)$ define non-$\dg$ functors.

\begin{corollary}
Let $f:X\hookrightarrow \PP^{n+1}$ be a smooth degree $d$-hypersurface. Then we have that the space of non-Fourier-Mukai functors arising as deformation of the shape $\psi_\eta$ can be computed explicitly  .
\end{corollary}

We now give a few examples of Twisted Hodge diamonds and the resulting dimensions for non-$\dg$ deformations of the closed immersions we computed using the Sage package \cite{TwistedHodgeDiamondsSage}. 

\begin{example}
Let $f:X\hookrightarrow \PP^{7}$ be a smooth degree $7$ hypersurface then its Hodge diamond is

$$\begin{matrix}
\;&\;&\;&\;&\;&\;&\;&1&\;&\;&\;&\;&\;\\
\;&\;&\;&\;&\;&\;&0&\;&0&\;&\;&\;&\;\\
\;&\;&\;&\;&\;&0&\;&1&\;&0&\;&\;&\;\\
\;&\;&\;&\;&0&\;&0&\;&0&\;&0&\;&\;\\
\;&\;&\;&0&\;&0&\;&1&\;&0&\;&0&\;\\
\;&\;&0&\;&0&\;&0&\;&0&\;&0&\;&0\\
\;&0&\;&\color{purple}\mathclap{36}\color{purple}&\;&\mathclap{2472}&\;&\mathclap{8093}&\;&\mathclap{2472}&\;&\mathclap{36}&\;&0\\
\;&\;&0&\;&0&\;&0&\;&0&\;&0&\;&0\\
\;&\;&\;&0&\;&0&\;&1&\;&0&\;&0&\;&\;&\;\\
\;&\;&\;&\;&0&\;&0&\;&0&\;&0&\;&\;&\;&\;\\
\;&\;&\;&\;&\;&0&\;&1&\;&0&\;&\;&\;&\;&\;\\
\;&\;&\;&\;&\;&\;&0&\;&0&\;&\;&\;&\;&\;&\;\\
\;&\;&\;&\;&\;&\;&\;&1.&\;&\;&\;&\;&\;&\;&\;\\
\end{matrix}$$

And so we have by \cite[Theorem~4.14]{Kueng2022},
\begin{align*}
\dim \ker\left(f_*:\HH^{10}\left(X,\OO_X\left(1\right)\right)\to \HH^{10}\left(\PP^7,f_*\OO_X\left(1\right)\right)\right)&=36\\
\end{align*}
All non-zero $\eta$ in this $36$-dimensional space give rise to non-Fourier-Mukai functors.
\end{example}

\begin{example}
Let $f:X\hookrightarrow \PP^{7}$ be a smooth degree $7$ hypersurface then its $8$ twisted Hodge diamond is

$$\begin{matrix}
\;&\;&\;&\;&\;&0&\;&\;&\;&\;&\;\\
\;&\;&\;&\;&0&\;&0&\;&\;&\;&\;\\
\;&\;&\;&0&\;&0&\;&0&\;&\;&\;\\
\;&\;&0&\;&0&\;&0&\;&0&\;&\;\\
\;&0&\;&0&\;&0&\;&0&\;&0&\;\\
\mathllap{299}6&\;&\color{purple}\mathclap{20993}\color{purple}&\;&\mathclap{15267}&\;&\mathclap{917}&\;&0&\;&0\\
\;&\mathllap{157}5&\;&0&\;&0&\;&0&\;&0&\;\\
\;&\;&\mathllap{577}5&\;&0&\;&0&\;&0&\;&\;\\
\;&\;&\;&\mathllap{1039}5&\;&0&\;&0&\;&\;&\;\\
\;&\;&\;&\;&\mathllap{900}2&\;&0&\;&\;&\;&\;\\
\;&\;&\;&\;&\;&\mathllap{299}6.&\;&\;&\;&\;&\;\\
\end{matrix}$$

And so we have by \cite[Theorem~4.14]{Kueng2022},
\begin{align*}
\dim \ker\left(f_*:\HH^{8}\left(X,\OO_X\left(-8\right)\right)\to \HH^{8}\left(X,f_*\OO_X\left(-8\right)\right)\right)&=\color{purple}20993\color{black}.\\
\end{align*}

Again every single $0\neq\eta$ in this $20993$-dimensional space give rise to a non-Fourier-Mukai functors.
\end{example}

Finally, paired with \cite[Proposition~4.16]{Kueng2022} we get the following corollary

\begin{corollary}
Let $f:X\hookrightarrow \PP^{n+1}$ be a smooth degree $d$-hypersurface. Then there exists a non-Fourier-Mukai functor
$$\psi_\eta: \D^b\left(\Qcoh X\right)\to \D^b\left(\Qcoh \PP^{n+1}\right).$$
\end{corollary}

Finally recall that all the above constructions are not liftable in the spectral setting as well by the following proposition which is the main reason we considered a field of characteristic zero.

\begin{proposition}\cite[Proposition~B.2.1.]{Rizzardo2019}\label{Proposition spectral lift}
If $\kk=\mathbb{Q}$ then the functor $\psi_\eta$ does not admit a spectral lift.
\end{proposition}

 \bibliographystyle{alpha}
\bibliography{Injectives_obstruct_Fourier-Mukai_functors}  

\end{document}

%% file: tikzstyledefs.tex
\tikzstyle heightone=[scale=.7,shift={(0,-.3)}]
\tikzstyle heightones=[scale=.8,xscale=.35,shift={(0,.1)}]
\tikzstyle heightoneonehalf=[scale=.9,shift={(0,-.2)}]
\tikzstyle heighttwo=[scale=.9,shift={(0,-.4)}]
\tikzstyle heighttwos=[scale=.5,xscale=.6,shift={(0,-.1)}]
\tikzstyle heightthree=[scale=.6,shift={(0,-.9)}]
\tikzstyle heightthrees=[scale=.4,xscale=.7,shift={(0,-.2)}]

\tikzstyle arrowstyle=[blue,semitransparent,scale=2]

\tikzstyle basiclabel=[draw=none,fill=none,shape=rectangle,inner sep=2pt,scale=.8]
\tikzstyle leftlabel=[basiclabel,anchor=east]
\tikzstyle rightlabel=[basiclabel,anchor=west]
\tikzstyle bottomlabel=[basiclabel,anchor=north]
\tikzstyle toplabel=[basiclabel,anchor=south]

\tikzstyle vertex=[circle,draw,fill=black,inner sep=1pt]
\tikzstyle ciliation=[circle,draw=none,fill=red,inner sep=1pt,semitransparent]
\tikzstyle ciliatednode=[vertex,pin={[pin distance=1mm,pin edge={semitransparent,red},ciliation]#1:{}}]

\tikzstyle vector=[black,thick,rectangle,draw=gray!50!yellow,top color=yellow!30,bottom color=black!10,scale=.8,inner sep=2pt]
\tikzstyle small vector=[vector,scale=.8]
\tikzstyle plain vector=[rectangle,draw=none,fill=white,scale=.7]

\tikzstyle my signal=[black,thick,signal,signal pointer angle=120,draw=blue!50,top color=blue!20,bottom color=black!10,scale=.8,inner sep=2pt]
\tikzstyle matrix=[my signal,signal from=south,signal to=north]
\tikzstyle reverse matrix=[my signal,signal from=north,signal to=south]

\tikzstyle small matrix=[matrix,scale=.7]
\tikzstyle reverse small matrix=[reverse matrix,scale=.7]
\tikzstyle matrix on edge=[small matrix,sloped,rotate=-90]
\tikzstyle reverse matrix on edge=[small matrix,sloped,rotate=90]

\tikzstyle trivalent=[very thick]
\tikzstyle dotdotdot=[decorate,decoration={markings,
    mark=at position .3 with{\node{.};},
    mark=at position .5 with {\node{.};},
    mark=at position .7 with {\node{.};}}]

\tikzstyle wavyup=[out=90,in=-90]
\tikzstyle wavydown=[out=-90,in=90]

\tikzstyle symmetrizer=[rectangle,fill=gray!10,draw=black]
\tikzstyle permutation=[symmetrizer]
\tikzstyle antisymmetrizer=[rectangle,fill=black,draw=black]
\tikzstyle symlabel=[draw=none,fill=none,black,scale=.8]
\tikzstyle asymlabel=[draw=none,fill=none,white,scale=.8]